\documentclass[preprint,12pt]{elsarticle}
\usepackage{amssymb}
\usepackage{amsthm}
\usepackage[figuresright]{rotating}
\usepackage{amsfonts}
\usepackage{amsmath}
\usepackage{latexsym}
\usepackage{euscript}
\usepackage{enumitem}
\usepackage{graphicx}
\usepackage{tikz}
\usepackage{subfigure}
\usepackage{caption}
\usepackage{float}
\usepackage{subfloat}
\usepackage{xcolor}
\newtheorem{theorem}{Theorem}[section]
\newtheorem{definition}{Definition}[section]
\newtheorem{lemma}{Lemma}[section]

\begin{document}
\begin{frontmatter}
\title{Mixed finite element method for a beam equation with the $p$-biharmonic operator}

\author[ubi]{Rui M. P. Almeida}
\ead{ralmeida@ubi.pt}
\author[ubi]{José C. M. Duque}
\ead{jduque@ubi.pt}
\author[uff]{Jorge Ferreira}
\ead{ferreirajorge2012@gmail.com}
\author[ubi]{Willian S. Panni\corref{cor1}}
\ead{willian.panni@ubi.pt}
\cortext[cor1]{Corresponding author.}

\address[ubi]{University of Beira
Interior, Mathematics and Applications Center, Rua Marquês d'Ávila e Bolama, 6201-001, Covilhã, Portugal.}

\address[uff]{Federal Fluminense
University, Department of Exact Sciences, Av. dos Trabalhadores 420,
27255-125, Volta Redonda, Brazil.}

\begin{abstract}
In this paper, we consider a nonlinear beam equation with the $p$-biharmonic operator, where $1 < p < \infty$. Using a change of variable, we transform the problem into a system of differential equations and prove the existence, uniqueness and regularity of the weak solution by applying the Lax-Milgram theorem and classical results of functional analysis. We investigate the discrete formulation for that system and, with the aid of the Brouwer theorem, we show that the problem has a discrete solution. The uniqueness and stability of the discrete solution are obtained through classical methods. After establishing the order of convergence, we apply the mixed finite element method to obtain an algebraic system of equations.  Finally, we implement the computational codes in Matlab software and perform the comparison between theory and simulations.
\end{abstract}

\begin{keyword}
$p$-biharmonic operator \sep weak solution \sep convergence order \sep mixed finite element method \sep numerical simulations \sep beam equation.
\hfill
\MSC[2020] 35A01 \sep 35D30 \sep 35J40 \sep 65N30.
\end{keyword}
\end{frontmatter}

\newpage
\section{Introduction}\label{sec1}

Let $\Omega$ be a bounded domain in $\mathbb{R}^N$ $(N \geq 1)$ with smooth boundary $\partial \Omega$. We consider the problem
\begin{equation}
\left\{
\begin{array}{ll}
\Delta^{2}_{p} u = f\left(x\right), & \mbox{in } \Omega, 
\\
u = 0, \Delta u = 0, & \mbox{on } \partial \Omega,
\end{array}
\right.  \label{Problema}
\end{equation}
\noindent where $\Delta^{2}_{p}$ is the fourth-order operator, called the $p$-biharmonic operator, defined by 
\begin{equation*}
\Delta^{2}_{p} u = \Delta \left( | \Delta u |^{p-2} \Delta u \right),
\end{equation*}
$p\in \mathbb{R}$ satisfies $1<p<\infty $ and  $f \in L^{2}\left( \Omega \right)$.

Nonlinear differential equations whose structure depends on the solution itself can arise in mathematical modeling of various real-life processes. Most of these models fall into the class of evolutionary equations with non-standard growth. Recently, differential equations and variational problems with non-standard growth conditions have attracted increasing attention (see, for example, \cite{Antontsev_Shmarev_2015} for a review). These problems arise in various branches of applied mathematics and physics, such as electro-rheological or thermo-rheological fluid flows \cite{Ruzicka_2000}, elastic mechanics \cite{Zhikov_1986} and digital image processing \cite{Chen_Levine_Rao_2006}. Solutions to these types of problems may exhibit interesting properties, such as finite time extinction, explosion, propagation of perturbations with finite velocity or waiting time phenomena, which are intrinsic to the solution of other nonlinear problems.

Fourth order PDEs have several applications, for example, in elastic beam deformations \cite{Gyulov_Morosanu_2010, Lazer_McKenna_1990}, in thin plate theory \cite{Danet_2018, Danet_2019} and in image processing \cite{You_Kaveh_2000}. The $p$-biharmonic problems are at the intersection of these fields of study.

For $p=2$, we refer the reader to \cite{Ciarlet_Raviart_1974}, where Ciarlet and Raviart investigate the biharmonic problem
\begin{equation}\label{Eq_p=2}
\Delta^2 u = f.
\end{equation}
In particular, they explained the standard variational formulation for the continuous problem and demonstrated that the solution exists and is unique. They used a standard procedure in optimization theory to replace the problem with an unrestricted one, where it was necessary to find a Lagrange multiplier and the function $u$. Furthermore, they estimated the error and convergence order with $C^{0}$-finite-elements associated with polynomials of degree $r-1$.

Behrens and Guzm\'{a}n \cite{Behrens_Guzman_2011} introduced a new mixed method for Problem \eqref{Eq_p=2}, based on a formulation where the problem was rewritten as a system of four first-order equations. A hybrid form of the method was introduced, which permitted a reduction of the globally coupled degrees of freedom to only those associated with Lagrange multipliers which approximate the solution and its derivative at the faces of the triangulation.

Pryer \cite{Pryer_2014} studied discontinuous Galerkin approximations of Problem \eqref{Problema} for $1 < p < \infty$ from a variational perspective. The author proposed a discrete variational formulation of the problem based on an appropriate definition of a finite element Hessian and studied the convergence of the method (without rates) using a semicontinuity argument.

In \cite{Candito_Bisci_2012}, Candito and Bisci's approach was based mainly on critical point theory. They showed that there are multiple weak solutions for Problem \eqref{Problema}, with $f=\lambda g(x,u)$ and $\lambda$ a positive parameter.

Lu and Fu \cite{Lu_Fu_2020} studied Problem \eqref{Problema} with $f=\lambda|u|^{p-2}u-|u|^{q-2}$. They proved that the problem has infinitely many solutions when $p<q<p^{\ast }=\frac{Np}{N-2p}$ and obtained a multiplicity existence result when $q=p^{\ast}$.

Chaharlang and Razani \cite{Chaharlang_Razani_2019} considered Problem \eqref{Problema} with $f=\mu \frac{|u|^{p-2}u}{|x|^{2p}}+\lambda g(x,u)$, where $\mu $ and $\lambda $ are real parameters. Applying the Ricceri critical point principle, variational methods and the Rellich inequality, they established the existence of at least one weak solution in two different cases of the nonlinear term at the origin.

More recently, Katzourakis and Pryer \cite{Katzourakis_Pryer_2019} studied the problem
\begin{equation*}
\Delta _{p}^{2}u=0. 
\end{equation*}
For a fixed $p$, they proposed a method based on C$^{0}$-mixed finite elements. They rewrote the minimization problem in a mixed formulation and proved, with norms that depend on the mesh, that the method converges under the minimum regularity of the solution. Moreover, under additional regularity assumptions and using an inf-sup condition, they showed that the approximation converges with specific rates that depend on $p$. They proved the convergence of the numerical solution to the weak solution of $\infty$-bilaplacian problem. Finally, they provided numerical experiments to validate their theoretical results.

In the present paper, we generalize the results of Katzourakis and Pryer in \cite{Katzourakis_Pryer_2019}, in the sense that we use classical norms that do not depend on the mesh. In addition, we prove the order of convergence, which is found to be optimal for some values of $p$. The paper is organized as follows. In Section $\ref{sec2}$, we present some known results concerning Lebesgue and Sobolev spaces that shall be required. In Section $\ref{sec3}$, we prove the existence and uniqueness of the weak solution to Problem $(\ref{Problema})$. In Section $\ref{sec4}$, we study the discretized problem using the finite element method. In particular, we prove the existence, uniqueness and regularity of the discrete solution and study the order of convergence. In Section $\ref{sec5}$, we present numerical simulations using Matlab software. The conclusions of the paper are presented in Section $\ref{sec6}$.

\section{Preliminaries}\label{sec2}

In this work, we use the standard notation for the Lebesgue and Sobolev spaces, $L^{p}\left( \Omega \right) $ and $W^{m,p}\left( \Omega \right)$, respectively  (see \cite{Adams_Fournier_2003}). We denote the inner product in $L^{2}\left( \Omega \right) $ by $(\cdot \mbox{, }\cdot )$ and, throughout this work, constants which are independent of the parameters and the functions involved are denoted by $C$.

The closure of $C_{0}^{\infty }(\Omega )$ in $W^{m,p}(\Omega )$ is denoted by $W_{0}^{m,p}(\Omega )$. For case $p=2$, we let $W^{m,2}(\Omega)=H^{m}(\Omega )$ and $W_{0}^{m,2}(\Omega )=H_{0}^{m}(\Omega )$. An equivalent norm of $H_{0}^{1}(\Omega )$ is given by $\left\Vert \nabla u\right\Vert_{L^{2}(\Omega )}$.

For the sake of clarity, we first recall some known results.

\begin{theorem}\label{Teo:u_H2}\textup{\protect\cite{Evans_2010}.} Let $\Omega \subset \mathbb{R}^N$ be a bounded and open set with $\partial \Omega$ of class $C^2$. If $f \in L^{2}(\Omega)$ and $u \in H_{0}^{1}(\Omega)$ satisfies
\begin{equation*}
\int_{\Omega }\nabla u\nabla \varphi dx=\int_{\Omega }f\varphi dx,\forall \ \varphi \in H_{0}^{1}(\Omega ),
\end{equation*}
then $u\in H^{2}(\Omega )$ and
\begin{equation*}
\left\Vert u \right\Vert_{H^{2}(\Omega)} \leq C \left\Vert f \right\Vert_{L^{2}(\Omega)},
\end{equation*}
where $C$ is a constant depending only on $\Omega$.
\end{theorem}

\begin{theorem}\label{Teo_Brouwer} \textup{\protect\cite{Mawhin_2020}.} Let $\overline{B}_{\varrho }(0)$ be a closed ball with center at the origin and radius $\varrho $. If $f:\overline{B}_{\varrho}(0)\rightarrow \mathbb{R}^{N}$ is continuous and $\left( f(x),x\right) \geq 0$ for every $x\in \partial \overline{B}_{\varrho }(0)$, then $f$ has a zero in $\overline{B}_{\varrho }(0)$.
\end{theorem}

\begin{theorem}
\label{Teo_mu} Let $0 < \mu \leq 1$ and $\Omega \subset \mathbb{R}^N$. If $u \in L^{2}(\Omega)$, then $\left|u\right|^{\mu} \in L^{1}(\Omega)$ and
\begin{equation*} 
\int_{\Omega}\left|u(x)\right|^{\mu}dx \leq C \left\Vert u \right\Vert_{L^{2}\left(\Omega\right)}^\mu,
\end{equation*}
where $C = \left|\Omega\right|^{\frac{2-\mu}{2}}$.
\end{theorem}

\begin{proof}
We have using the H\"{o}lder inequality, 
\begin{equation*}
\int_{\Omega}\left|u(x)\right|^{\mu}dx = 		\int_{\Omega}1\left|u(x)\right|^{\mu}dx \leq \left\Vert 1 \right\Vert_{L^{q}(\Omega)}\left\Vert \left|u\right|^{\mu} \right\Vert_{L^{p}(\Omega)},
\end{equation*}
where $p$ and $q$ are conjugate exponents, that is, $q=\frac{p}{p-1}$.
	
Since $2 \leq \frac{2}{\mu} < \infty$, we choose $p = \frac{2}{\mu}$, which implies $q = \frac{2}{2 - \mu}$. Substituting these values of $p$ and $q$ in the inequality above, we obtain
\begin{equation*}
\int_{\Omega}\left|u(x)\right|^{\mu}dx \leq \left\Vert 1 \right\Vert_{L^{\frac{2}{2 - \mu}}(\Omega)}\left\Vert \left|u\right|^{\mu} \right\Vert_{L^{\frac{2}{\mu}}(\Omega)}.
\end{equation*}
	
Noticing that
\begin{equation*}
\left\Vert 1 \right\Vert_{L^{\frac{2}{2 - \mu}}(\Omega)} =  \left(\int_{\Omega}\left|1\right|^{\frac{2}{2 - \mu}}dx\right)^{\frac{2 - \mu}{2}} = \left(\int_{\Omega} 1 dx\right)^{\frac{2 - \mu}{2}} = \left|\Omega\right|^{\frac{2-\mu}{2}} = C
\end{equation*}
and
\begin{equation*}
\left\Vert \left|u\right|^{\mu} \right\Vert_{L^{\frac{2}{\mu}}(\Omega)} =  \left(\int_{\Omega}\left|\left|u(x)\right|^{\mu}\right|^{\frac{2}{\mu}}dx\right)^{\frac{\mu}{2}} = \left(\left(\int_{\Omega} \left|u(x)\right|^{2} dx\right)^{\frac{1}{2}}\right)^{\mu} = \left\Vert u \right\Vert_{L^{2}\left(\Omega\right)}^\mu,
\end{equation*}
we have the required result.	
\end{proof}

We denote by $\mathcal{T}_{h}$ a non-degenerate partition of the polygonal domain $\Omega \subset \mathbb{R}^{N}$ in simplexes with parameter $h$, that is, the set $\Omega $ is subdivided into a finite number of subsets $T_{k}$, $k=1,\ldots ,n$, called finite elements, such that the following conditions are met:
\begin{itemize}
\item[i)] $\Omega = \displaystyle\cup_{k = 1}^{n} T_k$;

\item[ii)] $int(T_k) \neq \emptyset, \forall \ T_k \in \mathcal{T}_{h}$;

\item[iii)] $int(T_i) \cap int(T_j) = \emptyset, \forall \ T_i, T_j \in \mathcal{T}_{h}$ with $i \neq j$;

\item[iv)] Each side of $T_{k}$ either belongs to the boundary of the domain $\Omega $ or is a side of another $T_{i}\in \mathcal{T}_{h}$;

\item[v)] Each $T_k$ has a Lipschitz-continuous boundary.
\end{itemize}

We denote by $V_{h}$ the space of continuous functions on the closure $\overline{\Omega }$ of $\Omega $, which are polynomials of degree $r-1$, with $r\geq 2$, in each interval of $\mathcal{T}_{h}$ and vanish on $\partial \Omega $, that is,
\begin{equation*}
V_{h}=\left\{ u\in C_{0}^{0}\left( \overline{\Omega }\right) ;u_{\left\vert_{T_{k}}\right. }\mbox{is a polynomial of degree }r-1\mbox{ for all } T_{k}\in \mathcal{T}_{h}\right\}.
\end{equation*}

\begin{theorem}\label{Teo.Est.Interpolador} \textup{\protect\cite{Thomee_2006}.} If $V_{h} \subset H_{0}^{1}\left(\Omega\right)$ and $I_h: H^{r}(\Omega) \cap H_{0}^{1}\left(\Omega\right) \rightarrow V_h$ is the interpolation operator, then
\begin{equation*}
\left\Vert u - I_{h}u \right\Vert_{L^{2}(\Omega)} + h\left\Vert \nabla \left(u - I_{h}u\right)\right\Vert_{L^{2}(\Omega)} \leq Ch^{s} \left\Vert u \right\Vert_{H^{s}(\Omega)}, \ 1 \leq s \leq r,
\end{equation*}
where $u \in H^{s}(\Omega) \cap H_{0}^{1}\left(\Omega\right)$ and $C$ is a constant.
\end{theorem}

\begin{theorem}\label{Teo_Estabilidade} \textup{\protect\cite{Brenner_Scott_2008}.} If $u\in W^{m,p}(\Omega )$, with $0\leq m\leq r$ and $1\leq p\leq \infty $, then
\begin{equation*}
\left(\sum_{k=1}^{n} \left\Vert I_h u\right\Vert_{W^{m,p}(T_k)}^{p} \right)^{\frac{1}{p}} \leq C \left(\sum_{\left|\alpha\right| = m} \left\Vert D^{\alpha} u\right\Vert_{L^{p}(\Omega)}^{p} \right)^{\frac{1}{p}}.
\end{equation*}

In particular,
\begin{equation*}
\sum_{k=1}^{n} \left\Vert I_h u\right\Vert_{W^{m,p}(T_k)}^{p} \leq C \left\Vert u\right\Vert_{W^{m,p}(\Omega)}^{p} .
\end{equation*}
\end{theorem}

\begin{theorem}\label{Teo_Des} \textup{\protect\cite{Barrett_Liu_1993}.} For all $q > 1$ and $\delta \geq 0$, there is positive constant $C$ such that, for all $\xi, \kappa \in \mathbb{R}$ with $\xi \neq \kappa$,
\begin{equation*}
\left|\left| \xi \right|^{q-2}\xi - \left|\kappa\right|^{q-2}\kappa \right| \leq C\left|\xi - \kappa\right|^{1-\delta}\left(\left|\xi\right| + \left|\kappa\right| \right)^{q-2+\delta}.
\end{equation*}
\end{theorem}

\section{Existence and uniqueness of weak solution}\label{sec3}

In this section, we establish the existence and uniqueness of the weak solution to Problem \eqref{Problema}.

Following Katzourakis and Pryer \cite{Katzourakis_Pryer_2019}, we define the auxiliary variable
\begin{equation}  \label{v}
v = \left|\Delta u \right|^{p-2}\Delta u,
\end{equation}
whence
\begin{equation}
\left\vert v\right\vert ^{q-2}v=\Delta u.  \label{Delta_u}
\end{equation}

Using Equations \eqref{v} and \eqref{Delta_u}, we may rewrite Problem \eqref{Problema} as the following system of differential equations
\begin{equation}
\left\{
\begin{array}{ll}
\Delta v=f, & \mbox{ in }\Omega , 
\\
\Delta u=\left\vert v\right\vert ^{q-2}v, & \mbox{ in }\Omega , \\
u=0,\ v=0, & \mbox{ on }\partial \Omega .
\end{array}
\right.   \label{Problema_2}
\end{equation}

\begin{definition}
\label{Def_Sol_Fraca_p=cte} The pair $(u, v) \in H_{0}^{1}(\Omega) \times H_{0}^{1}(\Omega)$ is a weak solution to Problem \eqref{Problema_2} if, for all $(\psi, \eta) \in H_{0}^{1}(\Omega) \times H_{0}^{1}(\Omega)$, the following system is satisfied
\begin{equation}  \label{Formulacao_fraca}
\left\{
\begin{array}{ll}
\displaystyle\int_{\Omega} \nabla v \nabla \psi dx = -\int_{\Omega} f \psi dx, \vspace{0.1cm} &  
\\
\displaystyle\int_{\Omega} \nabla u \nabla \eta dx = - \int_{\Omega} \left|v\right|^{q-2}v\eta dx. &
\end{array}
\right.
\end{equation}
\end{definition}

The weak formulation \eqref{Formulacao_fraca} can be rewritten as
\begin{equation}  \label{Sitema_a}
\left\{
\begin{array}{ll}
a\left( v, \psi \right) = -\left(f, \psi\right), & \forall \ \psi \in H_{0}^{1}(\Omega), 
\\
a\left( u, \eta \right) = -\left(\left|v \right|^{q-2}v, \eta\right), & \forall \ \eta \in H_{0}^{1}(\Omega),
\end{array}
\right.
\end{equation}
where the bilinear form $a: H_{0}^{1}(\Omega) \times H_{0}^{1}(\Omega) \rightarrow \mathbb{R}$ is given by
\begin{equation}  \label{forma_bilinear_a}
a\left(u, v \right) = \int_{\Omega} \nabla u \nabla v dx, \ \forall \ u, v \in H_{0}^{1}(\Omega).
\end{equation}

Next, we use the Lax-Milgram theorem and the Sobolev embedding theorem with $\Omega \subset \mathbb{R}^N$  to proof the existence and uniqueness of the weak solution. Due to the embeddings, we first consider $N \leq 4$, and then $N > 4$.

\begin{theorem}\label{Teo_N<=4} Let $\Omega \subset \mathbb{R}^{N}$, $N \leq 4$, be a bounded and open set with $\partial \Omega $ of class $C^{2}$, $f\in L^{2}(\Omega )$ and $1<p<\infty $. Then there is a unique pair $(u,v)$ that is the weak solution to Problem \eqref{Problema_2}, in the sense of Definition \ref{Def_Sol_Fraca_p=cte}, moreover $(u,v)\in H^{2}(\Omega )\times H^{2}(\Omega )$.
\end{theorem}

\begin{proof}
Suppose $f \in L^{2}\left(\Omega\right) \subset H^{-1}\left(\Omega\right)$. Since the bilinear form $a$ defined in \eqref{forma_bilinear_a} is continuous and coercive, we have, by the Lax-Milgram theorem, that there is a unique $v \in H_{0}^1\left(\Omega\right)$ such that the first equation of System \eqref{Sitema_a} is satisfied.  According to Theorem \ref{Teo:u_H2}, $v \in H^{2}(\Omega)$. 
	
Now, we investigate $\left\vert v \right\vert^{q-2}v$, consider different cases.

If $\Omega \subset \mathbb{R}^N$, with $N < 4$ and $1 < p < \infty$, then $H^{2}(\Omega) \subset C ( \overline{\Omega})$ (see \cite[pp. $44$]{Novotny_Straskraba_2004}), whence $\left|v\right|^{q-2}v \in L^{\infty}(\Omega) \subset L^{2}(\Omega)$. 

If $\Omega \subset \mathbb{R}^4$ and $1 < p < 3$, then $1 < 2q-2 < \infty$. Thus,  $v \in H^{2}(\Omega)\subset L^{2q-2}(\Omega)$ (see \cite[pp. $44$]{Novotny_Straskraba_2004}), which implies  $\left|v\right|^{q-2}v \in L^{2}(\Omega)$. In fact,
\begin{equation}\label{Remark_1<p<=3}
\int_{\Omega} \left\vert\left|v\right|^{q-2}v\right\vert^{2}dx = \int_{\Omega}\left|v\right|^{2q-2}dx  \ \leq C,
\end{equation}
where $C$ is a constant.

If $\Omega \subset \mathbb{R}^4$ and $3 \leq p < \infty$, then $0 < 2q-2 \leq 1$. Since $v \in H^{2}(\Omega)\subset L^{2}(\Omega)$, it follows from Theorem \ref{Teo_mu} that
\begin{equation}\label{Remark_3<=p<infty}
\int_{\Omega}\left|v\right|^{2q-2}dx \leq C\left\Vert v \right\Vert_{L^{2}(\Omega)}^{2q-2},
\end{equation}
where $C$ is a constant. Thus, $\left|v\right|^{q-2}v \in L^{2}(\Omega)$.

Hence, we have  $\left|v\right|^{q-2}v \in L^{2}(\Omega)$ for $\Omega \subset \mathbb{R}^N$, with $N \leq 4$, and $1 < p < \infty$. As $L^{2}\left(\Omega\right) \subset H^{-1}\left(\Omega\right)$, we may conclude from the Lax-Milgram theorem that, for all $\left|v \right|^{q-2}v \in L^{2}\left(\Omega\right)$, there is a unique $u \in H_{0}^1\left(\Omega\right)$ such that the second equation of System \eqref{Sitema_a} is satisfied. Finally, using Theorem \ref{Teo:u_H2}, it follows that $u \in H^{2}(\Omega)$, as desired.
\end{proof}

\begin{theorem}\label{Teo_N>4} Let $\Omega \subset \mathbb{R}^{N}$, $N>4$, be a bounded and open set with $\partial \Omega $ of class $C^{2}$, $f\in L^{2}(\Omega )$ and $\frac{2N-4}{N}<p<\infty $. Then there is a unique pair $(u,v)$ that is the weak solution to Problem \eqref{Problema_2}, in the sense of Definition \ref{Def_Sol_Fraca_p=cte}, moreover $(u,v)\in H^{2}(\Omega )\times H^{2}(\Omega )$.
\end{theorem}

\begin{proof} The proof of the first equation is analogous to that of Theorem \ref{Teo_N<=4} and we have $v \in H^{2}(\Omega)$. Next, we investigate $\left|v\right|^{q-2}v$  considering different cases, as in Theorem \ref{Teo_N<=4}.
	
If $\frac{2N-4}{N} < p < 3$, then $1 < 2q-2 < \frac{2N}{N-4}$. Thus, $v \in H^{2}(\Omega)\subset L^{2q-2}(\Omega)$ (see \cite[pp. $44$]{Novotny_Straskraba_2004}) and, following the proof of the inequality in  \eqref{Remark_1<p<=3}, we have $\left|v\right|^{q-2}v \in L^{2}(\Omega)$.

If $3 \leq p < \infty$, then $0 < 2q-2 \leq 1$. We have $v \in H^{2}(\Omega)\subset L^{2}(\Omega)$, so reasoning as in the proof of Inequality \eqref{Remark_3<=p<infty}, we obtain $\left|v\right|^{q-2}v \in L^{2}(\Omega)$. 

Thus, for $\frac{2N-4}{N} < p < \infty$, we obtain $\left|v\right|^{q-2}v \in L^{2}(\Omega)$. Now, arguing as in the proof of Theorem \ref{Teo_N<=4}, we conclude that for the second equation of System \eqref{Sitema_a} there is a unique $u \in H_{0}^1\left(\Omega\right) \cap H^2\left(\Omega\right)$.
\end{proof}

\section{Discrete problem}\label{sec4}

In this section, we study the discrete problem associated with Problem \eqref{Problema_2}. For the sake of simplicity, we consider $\Omega \subset \mathbb{R}^{N}$ with $N=1$. For dimensions $N>1$, the proofs are similar, but rely on other Sobolev embeddings and require some restrictions on $p$.

\subsection{Existence, uniqueness and regularity of discrete solution}

We start by defining the concept of discrete solution to Problem \eqref{Problema_2}. Let $\mathcal{T}_{h}$ be a non-degenerate partition of the domain $\Omega $ with parameter $h$, and $V_{h}\subset H_{0}^{1}\left(\Omega \right) $ be the finite element space associated with $\mathcal{T}_{h}$ and formed by polynomials of degree at most $r-1$, with $r\geq 2$.

\begin{definition} \label{Def_Sol_Discreta_p=cte} The pair $(u_{h},v_{h})\in V_{h}\times V_{h}$ is said to be a discrete solution to Problem \eqref{Problema_2} if,  for every pair $(\psi _{h},\eta _{h})\in V_{h}\times V_{h}$, the following system is satisfied
\begin{equation}  \label{Sistema_discreto}
\left\{
\begin{array}{ll}
\left(\nabla v_h, \nabla \psi_h \right) = -\left(f, \psi_h \right), & \forall \ \psi_h \in V_{h}, 
\\
\left( \nabla u_h, \nabla\eta_h \right) = -\left(\left|v_h \right|^{q-2}v_h, \eta_h\right), & \forall \ \eta_h \in V_{h}.
\end{array}
\right.
\end{equation}
\end{definition}

\begin{lemma}
\label{Lemma_vh} If $v_{h}\in H_{0}^{1}(\Omega )$, then $\left\vert v_{h}\right\vert ^{q-1}\in L^{2}(\Omega )$ for $q>1$ and
\begin{equation}  \label{Eq_vh^q-1}
\left\Vert \left|v_h\right|^{q-1} \right\Vert_{L^{2}(\Omega)} \leq C
\left\Vert v_h \right\Vert_{H_{0}^{1}(\Omega)}^{q-1}.
\end{equation}
\end{lemma}

\begin{proof} By hypothesis, $v_h \in H_{0}^{1}(\Omega)$ and we know that $H_{0}^{1}(\Omega) \subset L^{\infty}(\Omega)  \subset L^{2}(\Omega)$. Thus
\begin{equation}\label{Txt01}
\left\Vert v_h \right\Vert_{L^{2}(\Omega)} \leq C \left\Vert v_h \right\Vert_{L^{\infty}(\Omega)}.
\end{equation}
On the other hand,
\begin{equation}\label{Txt02}
\left\Vert \left| v_h\right|^{q-1} \right\Vert_{L^{2}(\Omega)} = \left(\int_{\Omega} \left| \left| v_h\right|^{q-1}\right|^{2} dx\right)^{\frac{1}{2}} = \left( \int_{\Omega} \left| v_h\right|^{2q-2} dx \right)^{\frac{1}{2}}.
\end{equation}
From $(\ref{Txt01})$ and $(\ref{Txt02})$,
\begin{equation*}
\left\Vert \left| v_h\right|^{q-1} \right\Vert_{L^{2}(\Omega)} \leq \left(\left\Vert v_h \right\Vert_{L^{\infty}(\Omega)}^{2q-2} \right)^{\frac{1}{2}} \left( \int_{\Omega}1 dx \right)^{\frac{1}{2}} = C \left\Vert v_h \right\Vert_{L^{\infty}(\Omega)}^{q-1},
\end{equation*}
where $C = \left\vert \Omega \right\vert^{\frac{1}{2}} $. From the Sobolev embedding theorem, we now obtain \eqref{Eq_vh^q-1}. \end{proof}

Next, we show that Problem \eqref{Problema_2} has a discrete solution.

\begin{theorem}\label{Teo_existencia_uh_vh} If $f\in L^{2}(\Omega )$, then there is a discrete solution $(u_{h},v_{h})\in V_{h}\times V_{h}$ to Problem \eqref{Problema_2}, in the sense of Definition \ref{Def_Sol_Discreta_p=cte}.
\end{theorem}

\begin{proof} We define a continuous linear map $f_1: V_h \rightarrow V_h$ by
\begin{equation}\label{Eq:f1}
\left(f_1(v_h), \psi_h\right) = \left(\nabla v_h, \nabla \psi_h\right) + \left(f, \psi_h\right), \ \forall \ \psi_h \in V_h.
\end{equation}
	
Considering $\psi_h = v_h$, Equation \eqref{Eq:f1} becomes
\begin{equation*}
\left(f_1(v_h), v_h\right) = \left\Vert \nabla v_h \right\Vert_{L^{2}\left(\Omega\right)}^2 + \left(f, v_h\right).
\end{equation*}
	
Using the Poincaré and Young inequalities,
\begin{equation*}
\left(f_1(v_h), v_h\right) \geq \left(\dfrac{1}{C} - \varepsilon \right)\left\Vert v_h \right\Vert_{L^{2}\left(\Omega\right)}^2 - C_\varepsilon\left\Vert f \right\Vert_{L^{2}(\Omega)}^2,
\end{equation*}
where $C_\varepsilon = \frac{1}{4\varepsilon}$. In particular, to obtain  $\left(f_1(v_h), v_h\right) = 0$, we require that
\begin{equation*}
\left(\dfrac{1}{C} - \varepsilon \right)\left\Vert v_h \right\Vert_{L^{2}\left(\Omega\right)}^2 = C_\varepsilon\left\Vert f \right\Vert_{L^{2}(\Omega)}^2.
\end{equation*}
Choosing $\varepsilon$ such that $\left(\frac{1}{C} - \varepsilon \right) > 0$,
\begin{equation*}
\left\Vert v_h \right\Vert_{L^{2}\left(\Omega\right)} = \varrho_1\left\Vert f \right\Vert_{L^{2}(\Omega)},
\end{equation*}
where $\varrho_1 = \sqrt{\frac{C_\varepsilon C}{1-\varepsilon C}}.$
	
We define the closed ball $\overline{B}_{\varrho_2}(0) = \left\{ v_h \in V_h; \left\Vert v_h \right\Vert_{L^{2}\left(\Omega\right)} \leq   \varrho_2\right\}$, with $\varrho_2 = \varrho_1\left\Vert f \right\Vert_{L^{2}(\Omega)}$. So, for all $v_h \in \partial\overline{B}_{\varrho_2}(0)$, we have $\left(f_1(v_h), v_h\right) \geq 0$. Hence, by Theorem \ref{Teo_Brouwer}, there is $v_{h}^{*} \in \overline{B}_{\varrho_2}(0)$ such that $f_1(v_{h }^{*}) = 0$. Therefore,
\begin{equation*}
\left(\nabla v_{h}^{*}, \nabla \psi_h\right) = - \left(f, \psi_h\right), \ \forall \ \psi_h \in V_h,
\end{equation*}
whence $v_{h}^{*}$  is a solution of the first equation of Problem \eqref{Sistema_discreto}.
	
Analogously, we define a linear map $f_2: V_h \rightarrow V_h$ by
\begin{equation}\label{Eq:f2}
\left(f_2(u_h), \eta_h\right) = \left(\nabla u_h, \nabla \eta_h\right) + \left(\left\vert v_h\right\vert^{q-2}v_h, \eta_h\right), \  \forall \ \eta_h \in V_h,
\end{equation}
which is continuous. In fact, by the Minkowski inequality,
\begin{equation}\label{Eq:f2_modulo}
\left|\left(f_2(u_h), \eta_h\right)\right| \leq  \left|\left(\nabla u_h, \nabla \eta_h\right)\right| + \left| \left(\left\vert v_h\right\vert^{q-2}v_h, \eta_h\right)\right|.
\end{equation}
	
Using the H\"{o}lder inequality,
\begin{align}
\left|\left(\nabla u_h, \nabla \eta_h\right)\right|
& \leq \left(\int_{\Omega} \left\vert \nabla u_h\right\vert^2 dx\right)^\frac{1}{2} \left(\int_{\Omega} \left\vert \nabla \eta_h\right\vert^2 dx\right)^\frac{1}{2}
\nonumber \\
& = \left\Vert \nabla u_h\right\Vert_{L^{2}(\Omega)} \left\Vert \nabla \eta_h\right\Vert_{L^{2}(\Omega)}
\nonumber \\
& \leq C. \label{Eq:f2_c1}
\end{align}
	
By the embeddings $H_{0}^{1}(\Omega) \subset L^{\infty}(\Omega) \subset L^{2}(\Omega)$ and Lemma \ref{Lemma_vh}, 
\begin{align}
\left|\left(\left\vert v_h\right\vert^{q-2}v_h, \eta_h\right)\right| & \leq \int_{\Omega} \left\vert \left\vert v_h\right\vert^{q-2}v_h\right\vert \left\vert\eta_h  \right\vert dx 
\nonumber \\
& = \int_{\Omega} \left\vert \left\vert v_h\right\vert^{q-1}\right\vert \left\vert\eta_h \right\vert dx 
\nonumber \\
& \leq \left\Vert \left\vert v_h\right\vert^{q-1} \right\Vert_{L^{2}(\Omega)} \left\Vert \eta_h\right\Vert_{L^{2}(\Omega)}
\nonumber \\
& \leq C. \label{Eq:f2_c2}
\end{align}
	
Now, substituting \eqref{Eq:f2_c1} and \eqref{Eq:f2_c2} in \eqref{Eq:f2_modulo}, it follows that the linear map $f_2$ is continuous.
	
On the other hand, taking $\eta_h = u_h$ in Equation \eqref{Eq:f2},
\begin{equation*}
\left(f_2(u_h), u_h\right) = \left\Vert\nabla u_h\right\Vert_{L^{2}\left(\Omega\right)}^2 + \left(\left\vert v_h\right\vert^{q-2}v_h, u_h\right).
\end{equation*}
	
Using the Poincaré and Young inequalities,
\begin{equation*}
\left(f_2(u_h), u_h\right) \geq \left(\dfrac{1}{C} - \varepsilon \right)\left\Vert u_h \right\Vert_{L^{2}\left(\Omega\right)}^2 - C_\varepsilon\left\Vert \left\vert v_h \right\vert^{q-1} \right\Vert_{L^{2}(\Omega)}^2,
\end{equation*}
where $C_\varepsilon = \frac{1}{4\varepsilon}$. In particular, we require that $\left(f_2(v_h), u_h\right) = 0$, which implies 	
\begin{equation*}
\left(\dfrac{1}{C} - \varepsilon \right)\left\Vert u_h \right\Vert_{L^{2}\left(\Omega\right)}^2 = C_\varepsilon\left\Vert\left\vert v_h \right\vert^{q-1} \right\Vert_{L^{2}(\Omega)}^2.
\end{equation*}
	
If we choose $\varepsilon$ such that $\left(\frac{1}{C} - \varepsilon \right) > 0$ and define $\varrho_3 = \sqrt{\frac{C_\varepsilon C}{1-\varepsilon C}}$, then
\begin{equation*}
\left\Vert u_h \right\Vert_{L^{2}\left(\Omega\right)} = \varrho_3\left\Vert \left\vert v_h \right\vert^{q-1} \right\Vert_{L^{2}(\Omega)}.
\end{equation*}
	
Let $\overline{B}_{\varrho_4}(0)$ be the closed ball $\overline{B}_{\varrho_4}(0) = \left\{ u_h \in V_h; \left\Vert u_h \right\Vert_{L^{2}\left(\Omega\right)} \leq  \varrho_4\right\}$, with $\varrho_4 = \varrho_3\left\Vert \left\vert v_h \right\vert^{q-1} \right\Vert_{L^{2}(\Omega)}$. Then, for all $u_h \in \partial\overline{B}_{\varrho_4}(0)$, we obtain $\left(f_2(u_h), u_h\right) \geq 0$. So, by Theorem \ref{Teo_Brouwer}, there is a solution $u_{h}^{*} \in \overline{B}_{\varrho_4}(0)$ to the second equation of Problem \eqref{Sistema_discreto}.
\end{proof}

In the next theorem, we prove that the discrete solution to Problem \eqref{Problema_2} is unique.

\begin{theorem}
\label{Teo_unicidade_uh_vh} If $f\in L^{2}(\Omega )$, then the discrete solution $(u_{h},v_{h})\in V_{h}\times V_{h}$ to Problem \eqref{Problema_2}, in the sense of Definition \ref{Def_Sol_Discreta_p=cte}, is unique.
\end{theorem}

\begin{proof} Suppose $v_1, v_2 \in V_h$ are solutions to the first equation of Problem \eqref{Sistema_discreto}. Then
\begin{equation*}
\left(\nabla v_1, \nabla \psi_h\right) = - \left(f, \psi_h\right), \ \forall \ \psi_h \in V_h
\end{equation*}
and
\begin{equation*}
\left(\nabla v_2, \nabla \psi_h\right) = - \left(f, \psi_h\right), \ \forall \ \psi_h \in V_h.
\end{equation*}
	
Subtracting the previous equations and applying the inner product and gradient properties, we obtain
\begin{equation*}
\left(\nabla \left(v_1 - v_2\right), \nabla \psi_h\right) = 0.
\end{equation*}
	
The above equation is satisfied for all $\psi_h \in V_h$, so we can assume $\psi_h = v_1 - v_2$. Then
\begin{equation*}
\left\Vert \nabla \left(v_1 - v_2\right) \right\Vert_{L^{2}(\Omega)}^{2} = 0.
\end{equation*}
	
Using the Poincaré inequality,
\begin{equation*}
\left\Vert v_1 - v_2 \right\Vert_{L^{2}(\Omega)} = 0.
\end{equation*}
	
Consequently, $v_1 = v_2$ in $H_{0}^{1}(\Omega)$ and so there is a unique solution to the first equation of Problem \eqref{Sistema_discreto}.
	
Similarly, if $ u_1, u_2 \in V_h$ are solutions to the second equation of Problem \eqref{Sistema_discreto}, then
\begin{equation*}
\left(\nabla u_1, \nabla \eta_h\right) = - \left(\left|v_h\right|^{q-2}v_h, \eta_h \right), \ \forall \ \eta_h \in V_h
\end{equation*}
and
\begin{equation*}
\left(\nabla u_2, \nabla \eta_h\right) = - \left(\left|v_h\right|^{q-2}v_h, \eta_h \right), \ \forall \ \eta_h \in V_h.
\end{equation*}
	
As before, if we subtract the previous equations and apply the properties of the inner product and gradient, then  $\left(\nabla \left(u_1 - u_2\right), \nabla \eta_h\right) = 0.$ Since the above equation must be satisfied for all $\eta_h \in V_h$, we assume that $\eta_h = u_1 - u_2$, which implies $\left\Vert \nabla \left(u_1 - u_2\right) \right\Vert_{L^{2}(\Omega)}^{2} = 0.$ Thus, proceeding in a similar way as before, we conclude that $u_1 = u_2$ in $H_0^1(\Omega)$. This means that there is a unique solution to the second equation of Problem \eqref{Sistema_discreto}.
\end{proof}

In the next theorem, we prove the stability of the discrete solution $(u_{h},v_{h})$ to Problem \eqref{Problema_2}.

\begin{theorem}
\label{Teo_estabilidade_vh_uh} Let $(u_{h},v_{h})\in V_{h}\times V_{h}$ be the discrete solution of Problem \eqref{Problema_2}, in the sense of Definition \ref{Def_Sol_Discreta_p=cte}. Then, for all $f\in L^{2}(\Omega)$,
\begin{align}  \label{v_h < f}
\left\Vert v_h\right\Vert_{H_{0}^{1}(\Omega)} \leq C \left\Vert f
\right\Vert_{L^{2}(\Omega)}, 
\\
\left\Vert u_h\right\Vert_{H_{0}^{1}(\Omega)} \leq C \left\Vert f
\right\Vert_{L^{2}(\Omega)}^{q-1}.  \label{u_h < f}
\end{align}
\end{theorem}

\begin{proof} From the first equation of Problem \eqref{Sistema_discreto}, with $\psi_h = v_h$,
\begin{align*}
\left( \nabla v_h, \nabla v_h \right) = -\left( f, v_h \right) \Longrightarrow \left\Vert \nabla v_h \right\Vert_{L^{2}(\Omega)}^2 \leq \left\vert \left( f, v_h \right) \right\vert.
\end{align*}
Using the Young and Poincaré inequalities, 
\begin{equation*}
\left\Vert \nabla v_h \right\Vert_{L^{2}(\Omega)}^2 \leq C_\varepsilon \left\Vert f \right\Vert_{L^{2}(\Omega)}^2 + C\varepsilon\left\Vert \nabla v_h \right\Vert_{L^{2}(\Omega)}^2,
\end{equation*}
where $C_\varepsilon = \frac{1}{4\varepsilon}$. Choosing $\varepsilon$ such that $C\varepsilon = \frac{1}{2}$ and again, using the Poincaré inequality, we obtain \eqref{v_h < f}.
	
On the other hand, if $\eta_h = u_h$, then, from the second equation of Problem \eqref{Sistema_discreto},
\begin{equation*}
\left\Vert \nabla u_h \right\Vert_{L^{2}(\Omega)}^2 = -\left( \left| v_h \right|^{q-2} v_h, u_h \right),
\end{equation*}
which implies
\begin{equation*}
\left\Vert \nabla u_h \right\Vert_{L^{2}(\Omega)}^2 \leq \left\vert \left( \left| v_h \right|^{q-2} v_h, u_h \right) \right\vert \leq \int_{\Omega} \left| \left| v_h \right|^{q-1} \right| \left| u_h \right| dx.
\end{equation*}
By the Young and Poincaré inequalities,
\begin{equation*}
\left\Vert \nabla u_h \right\Vert_{L^{2}(\Omega)}^2 \leq C_\varepsilon \left\Vert  \left| v_h \right|^{q-1} \right\Vert_{L^{2}(\Omega)}^2 + C\varepsilon\left\Vert \nabla u_h \right\Vert_{L^{2}(\Omega)}^2,
\end{equation*}
where $C_\varepsilon = \frac{1}{4\varepsilon}$. Choosing $\varepsilon$ such that $C\varepsilon = \frac{1}{2}$, using the embedding $H_0^1(\Omega) \subset L^{\infty}(\Omega)$ and $(\ref{v_h < f})$, we obtain \eqref{u_h < f}.
\end{proof}

\subsection{Convergence order}

We now investigate the order of convergence for Problem \eqref{Sistema_discreto}. To this end, we first consider $1<p\leq 2$, which implies $2\leq q<\infty $, and then we consider $2<p<\infty $, which corresponds to $1<q<2$.

\begin{theorem}
\label{Teo_Ordem_Conv_1<p<2} Let $\Omega \subset \mathbb{R}^N$ be a bounded and open set with $\partial \Omega$ of class $C^2$ and  $1 < p \leq 2$. Under these conditions, for Problem \eqref{Sistema_discreto}, we have
\begin{eqnarray}
\left\Vert \nabla \left(v - v_h\right) \right\Vert_{L^2(\Omega)} & \leq & Ch^{s-1}\left\Vert v \right\Vert_{H^s(\Omega)},  \label{Ordem_grad v_1<p<2}
\\
\left\Vert v - v_h \right\Vert_{L^{2}(\Omega)} & \leq & Ch^s\left\Vert v \right\Vert_{H^{s}(\Omega)},  \label{Ordem_v_1<p<2} 
\\
\left\Vert \nabla\left(u-u_h\right)\right\Vert_{L^{2}(\Omega)} & \leq & Ch^{s-1}\left\Vert u \right\Vert_{H^{s}(\Omega)} + Ch^s\left\Vert v \right\Vert_{H^{s}(\Omega)},  \label{Ordem_grad u_1<p<2} 
\\
\left\Vert u - u_h \right\Vert_{L^{2}(\Omega)} & \leq & Ch^{s}\left\Vert u \right\Vert_{H^{s}(\Omega)} + C\left(h^{s} + h^{s+1}\right)\left\Vert v \right\Vert_{H^{s}(\Omega)},  \label{Ordem_u_1<p<2}
\end{eqnarray}
where $1 \leq s \leq r$ and $C$ is a constant.
\end{theorem}

\begin{proof} For the first equation of Problem \eqref{Sistema_discreto}, the proof is the same as that of Thomée \cite[Theorem $1.1$, pp. $5$]{Thomee_2006}, which gives \eqref{Ordem_grad v_1<p<2} and \eqref{Ordem_v_1<p<2}.
	
Using the same tools, subtracting the second equation of System \eqref{Sitema_a}, with $\eta = \eta_h \in V_h$, from the second equation of System \eqref{Sistema_discreto}, we obtain
\begin{equation}\label{s3}
\left(\nabla\left( u - u_h\right), \nabla \eta _h\right) = -(\left|v\right|^{q-2}v - \left|v_h\right|^{q-2}v_h, \eta_h).
\end{equation}
	
Let us now define
\begin{eqnarray}
\gamma = u - \tilde{u}_h,\label{gamma1}
\\
\theta = \tilde{u}_h - u_h,\label{theta1}
\end{eqnarray}
where $\tilde{u}_h$ represents the interpolation of $u$ in the space $V_h$. Thus,
\begin{equation*}
\left(\nabla \theta, \nabla \eta _h\right) = -\left(\nabla \gamma, \nabla \eta _h\right) -(\left|v\right|^{q-2}v - \left|v_h\right|^{q-2}v_h, \eta_h).
\end{equation*}
	
By definition, $\theta \in V_h$, so, considering $\eta_h = \theta$,
\begin{equation*}
\left\Vert \nabla \theta \right\Vert_{L^{2}(\Omega)}^2   \leq   \left\vert \left(\nabla \gamma, \nabla \theta\right) + (\left|v\right|^{q-2}v - \left|v_h\right|^{q-2}v_h, \theta) \right\vert.
\end{equation*}
	
Using the Minkowski, H\"{o}lder and Young inequalities, we obtain 	
\begin{align}
\left\Vert \nabla \theta \right\Vert_{L^{2}(\Omega)}^2 \leq & \  C_{\varepsilon_1}\left\Vert \nabla \gamma \right\Vert_{L^{2}(\Omega)}^2 + \varepsilon_1\left\Vert \nabla \theta \right\Vert_{L^{2}(\Omega)}^2
\nonumber \\
& + C_{\varepsilon_2}\left\Vert \left|v\right|^{q-2}v - \left|v_h\right|^{q-2}v_h \right\Vert_{L^{2}(\Omega)}^2 + \varepsilon_2\left\Vert \theta \right\Vert_{L^{2}(\Omega)}^2, \label{s4}
\end{align}
where $C_{\varepsilon_1} = \frac{1}{4\varepsilon_1}$ and $C_{\varepsilon_2} = \frac{1}{4\varepsilon_2}$. Using Theorem \ref{Teo_Des}, with $\delta = 0$, and the H\"{o}lder inequality, we obtain
\begin{equation*}
\left\Vert \left| v \right|^{q-2}v - \left|v_h\right|^{q-2}v_h \right\Vert_{L^{2}(\Omega)}^2\leq C\left\Vert v - v_h\right\Vert_{L^{2}(\Omega)}^2 \left\Vert \left(\left|v\right| + \left|v_h\right| \right)^{q-2}\right\Vert_{L^{2}(\Omega)}^2.
\end{equation*}
Since $v, v_h \in H_{0}^{1}(\Omega)$ and $\Omega \subset \mathbb{R}$, we have $H_{0}^{1}(\Omega) \subset L^{\infty}(\Omega)$. This implies $|v|$ and $|v_h|$ are bounded by a constant, whence
\begin{equation}\label{s5}
\left\Vert \left| v \right|^{q-2}v - \left|v_h\right|^{q-2}v_h \right\Vert_{L^{2}(\Omega)}^2 \leq C\left\Vert v - v_h\right\Vert_{L^{2}(\Omega)}^2.
\end{equation}
Substituting $(\ref{s5})$ in $(\ref{s4})$ gives
\begin{equation*}
\left\Vert \nabla \theta \right\Vert_{L^{2}(\Omega)}^2  \leq C_{\varepsilon_1}\left\Vert \nabla \gamma \right\Vert_{L^{2}(\Omega)}^2 + \varepsilon_1\left\Vert \nabla \theta \right\Vert_{L^{2}(\Omega)}^2 + C\left\Vert v - v_h\right\Vert_{L^{2}(\Omega)}^2 + \varepsilon_2\left\Vert \theta \right\Vert_{L^{2}(\Omega)}^2.
\end{equation*}
Applying the Poincaré inequality, choosing $\varepsilon_1$ and $\varepsilon_2$ such that $\varepsilon_1 + C\varepsilon_2 = \frac{1}{2}$, using Theorem \ref{Teo.Est.Interpolador} and Inequality $(\ref{Ordem_v_1<p<2})$, we obtain
\begin{equation}\label{d10}
\left\Vert \nabla \theta \right\Vert_{L^{2}(\Omega)}  \leq Ch^{s-1}\left\Vert u \right\Vert_{H^{s}(\Omega)}  + Ch^{s}\left\Vert v \right\Vert_{H^{s}(\Omega)}.
\end{equation}
	
Noting that
\begin{equation*}
\left\Vert  \nabla\left(u-u_h\right)\right\Vert_{L^{2}(\Omega)} = \left\Vert \nabla \gamma+ \nabla \theta \right\Vert_{L^{2}(\Omega)} \leq \left\Vert \nabla \gamma \right\Vert_{L^{2}(\Omega)} + \left\Vert \nabla \theta \right\Vert_{L^{2}(\Omega)}\label{23/07_Eq02}
\end{equation*}
and using Theorem \ref{Teo.Est.Interpolador} and \eqref{d10}, we conclude \eqref{Ordem_grad u_1<p<2}.
	
Finally, we determine the order of convergence for $u$. To this end, let  $\varphi \in L^{2}(\Omega)$ and $\psi \in H^{2}(\Omega) \cap H_{0}^{1} (\Omega)$ be the solution of $-\Delta \psi = \varphi$. Then,
\begin{equation*}
\left(u - u_h, \varphi\right) = - \left(u - u_h, \Delta \psi\right).
\end{equation*}
The next inequality follows from the Green formula, Equation \eqref{s3} and from the Minkowski and H\"{o}lder inequalities. 
\begin{align*}
\left(u - u_h, \varphi\right)  \leq &  \left\Vert \nabla\left(u - u_h\right)\right\Vert_{L^{2}(\Omega)}\left\Vert \nabla \left(\psi - \psi_h\right)\right\Vert_{L^{2}(\Omega)}
\\
& + \left\Vert \left|v\right|^{q-2}v - \left|v_h\right|^{q-2}v_h \right\Vert_{L^{2}(\Omega)}\left\Vert \psi_h\right\Vert_{L^{2}(\Omega)}.
\end{align*}
By $(\ref{Ordem_grad u_1<p<2})$, $(\ref{s5})$ and \eqref{Ordem_v_1<p<2}, 
\begin{align}
\left(u - u_h, \varphi\right) \leq &  \left(Ch^{s-1}\left\Vert u \right\Vert_{H^{s}(\Omega)}  + Ch^{s}\left\Vert v \right\Vert_{H^{s}(\Omega)}\right)\left\Vert \nabla \left(\psi - \psi_h\right)\right\Vert_{L^{2}(\Omega)}
\nonumber \\
& + Ch^s\left\Vert v \right\Vert_{H^{s}(\Omega)}\left\Vert \psi_h\right\Vert_{L^{2}(\Omega)}. \label{d14}
\end{align}
As $\psi_h \in V_h$ is arbitrary, we can consider $\psi_h$ as the interpolation of $\psi$ in the space $V_h$, that is, $\psi_h = I_h \psi = \tilde{\psi}_h$. Using Theorem \ref{Teo.Est.Interpolador} and the fact that $\psi \in H^{2}(\Omega)$, we obtain
\begin{equation}\label{d15}
\left\Vert \nabla \left(\psi - \tilde{\psi}_h\right)\right\Vert_{L^{2}(\Omega)}  \leq Ch \left\Vert \psi \right\Vert_{H^{2}(\Omega)}.
\end{equation}
Substituting $(\ref{d15})$ in $(\ref{d14})$, 
\begin{align}
\left(u - u_h, \varphi\right) \leq & \left(Ch^{s-1}\left\Vert u \right\Vert_{H^{s}(\Omega)}  + Ch^{s}\left\Vert v \right\Vert_{H^{s}(\Omega)}\right)Ch\left\Vert \psi \right\Vert_{H^{2}(\Omega)}
\nonumber \\
& +  Ch^s\left\Vert v \right\Vert_{H^{s}(\Omega)}\left\Vert \tilde{\psi}_h\right\Vert_{L^{2}(\Omega)}.  \label{d16}
\end{align}
Since $\psi \in H^{2}(\Omega)$, by Theorem \ref{Teo_Estabilidade}, $\left\Vert \tilde{\psi}_h \right\Vert_{H^{1}(\Omega)} \leq C \left\Vert \psi\right\Vert_{H^{2}(\Omega)}$. Furthermore, $H^{1}(\Omega) \subset L^{2}(\Omega)$, so
\begin{equation}
\left\Vert \tilde{\psi}_h\right\Vert_{L^{2}(\Omega)} \leq C \left\Vert \psi\right\Vert_{H^{2}(\Omega)}. \label{d20}
\end{equation}
Substituting \eqref{d20} in \eqref{d16}, we obtain the inequality 
\begin{align*}
\left(u - u_h, \varphi\right) \leq & \left(Ch^{s-1}\left\Vert u \right\Vert_{H^{s}(\Omega)}  + Ch^{s}\left\Vert v \right\Vert_{H^{s}(\Omega)}\right)Ch\left\Vert \psi \right\Vert_{H^{2}(\Omega)} 
\\
& +  Ch^s\left\Vert v \right\Vert_{H^{s}(\Omega)}C\left\Vert \psi\right\Vert_{H^{2}(\Omega)}.
\end{align*}
By Theorem \ref{Teo:u_H2}, $\left\Vert \psi \right\Vert_{H^{2}(\Omega)} \leq C \left\Vert \varphi \right\Vert_{L^{2}(\Omega)}$. Consequently,
\begin{equation*}
\left(u - u_h, \varphi\right) \leq C\left\Vert \varphi\right\Vert_{L^{2}(\Omega)}\left(h^{s}\left\Vert u \right\Vert_{H^{s}(\Omega)}  + \left(h^{s} + h^{s+1}\right)\left\Vert v \right\Vert_{H^{s}(\Omega)}\right).
\end{equation*}
Finally considering $\varphi = u - u_h$, we arrive at the inequality \eqref{Ordem_u_1<p<2}.
\end{proof}

We now consider $2<p<\infty $, so its conjugate exponent satisfies $1<q<2$.

\begin{theorem}
\label{Teo_Ordem_Conv_p>2} Let $\Omega \subset \mathbb{R}^N$ be a bounded and open set with $\partial \Omega$ of class $C^2$ and $2 < p < \infty$. For Problem \eqref{Sistema_discreto}, we have
\begin{align}
\left\Vert \nabla \left(v - v_h\right) \right\Vert_{L^2(\Omega)} & \leq Ch^{s-1}\left\Vert v \right\Vert_{H^s(\Omega)},  \label{Ordem_grad v_p>2} 
\\
\left\Vert v - v_h \right\Vert_{L^{2}(\Omega)} & \leq Ch^s\left\Vert v \right\Vert_{H^{s}(\Omega)},  \label{Ordem_v_p>2} 
\\
\left\Vert \nabla\left(u-u_h\right)\right\Vert_{L^{2}(\Omega)} & \leq Ch^{s-1}\left\Vert u \right\Vert_{H^{s}(\Omega)} + Ch^{\frac{s}{p-1}}\left\Vert v \right\Vert_{H^{s}(\Omega)}^{\frac{1}{p-1}},
\label{Ordem_grad u_p>2} 
\\
\left\Vert u - u_h \right\Vert_{L^{2}(\Omega)} & \leq Ch^s \left\Vert u\right\Vert_{H^{s}(\Omega)} + C\left(h^{\frac{s+p-1}{p-1}} + h^{\frac{s}{p-1}}\right) \left\Vert v \right\Vert_{H^{s}(\Omega)}^{\frac{1}{p-1}},
\label{Ordem_u_p>2}
\end{align}
where $1 \leq s \leq r$ and $C$ is a constant.
\end{theorem}

\begin{proof} The first equation of Problem \eqref{Sistema_discreto} does not depend on $p$ or $q$, so the proof is the same as that of Thomée \cite[Theorem $1.1$, pp. $5$]{Thomee_2006}. Hence we have \eqref{Ordem_grad v_p>2} and \eqref{Ordem_v_p>2}.
	
Now, we proceed in a manner analogous to Theorem \ref{Teo_Ordem_Conv_1<p<2}, until we obtain 
\begin{align}
\left\Vert \nabla \theta \right\Vert_{L^{2}(\Omega)}^2 \leq & \ C_{\varepsilon_1}\left\Vert \nabla \gamma \right\Vert_{L^{2}(\Omega)}^2 + \varepsilon_1\left\Vert \nabla \theta \right\Vert_{L^{2}(\Omega)}^2
\nonumber \\
& + C_{\varepsilon_2}\left\Vert \left|v\right|^{q-2}v - \left|v_h\right|^{q-2}v_h \right\Vert_{L^{2}(\Omega)}^2 + \varepsilon_2\left\Vert \theta \right\Vert_{L^{2}(\Omega)}^2. \label{ss4}
\end{align}
Then, applying Theorem \ref{Teo_Des} with $\xi = v$, $\kappa = v_h$, $\delta = \frac{p-2}{p-1} > 0$ and taking into account that $|v|$ and $|v_h|$ are bounded by a constant, we arrive at the inequality
\begin{equation}\label{sS5}
\left\Vert \left| v \right|^{q-2}v - \left|v_h\right|^{q-2}v_h \right\Vert_{L^{2}(\Omega)}^2 \leq C\left\Vert \left|v - v_h\right|^{\frac{1}{p-1}}\right\Vert_{L^{2}(\Omega)}^2.
\end{equation}
Using \eqref{sS5} in \eqref{ss4}, and applying the Poincaré inequality, choosing $\varepsilon_1$ and $\varepsilon_2$ so that $\varepsilon_1 + C\varepsilon_2 = \frac{1}{2}$, and taking in account Theorem \ref{Teo.Est.Interpolador}, we deduce
\begin{equation}\label{Eq04:07_06}
\dfrac{1}{2}\left\Vert \nabla \theta \right\Vert_{L^{2}(\Omega)}^2  \leq Ch^{2(s-1)}\left\Vert u \right\Vert_{H^{s}(\Omega)}^2  + C\left\Vert \left|v - v_h\right|^{\frac{1}{p-1}}\right\Vert_{L^{2}(\Omega)}^2.
\end{equation}
Note that
\begin{equation}
\left\Vert \left|v - v_h\right|^{\frac{1}{p-1}}\right\Vert_{L^{2}(\Omega)}^2  =  \int_{\Omega} \left|\left|v - v_h\right|^{\frac{1}{p-1}} \right|^2 dx =  \int_{\Omega} \left|v - v_h\right|^{\frac{2}{p-1}} dx. \label{Eq01:07_06}
\end{equation}
In order to analyze an estimate for Equation $(\ref{Eq01:07_06})$, it is necessary to consider $2 < p < 3$ and $p \geq 3$.
	
If $2 < p < 3$, then $1 < \frac{2}{p-1} < 2$. Consequently, we can define a norm in  $L^{\frac{2}{p-1}}\left(\Omega\right)$ and Equation $(\ref{Eq01:07_06})$ becomes
\begin{equation} \label{Eq03:07_06}
\left\Vert \left|v - v_h\right|^{\frac{1}{p-1}}\right\Vert_{L^{2}(\Omega)}^2 = \left\Vert v - v_h \right\Vert_{L^{\frac{2}{p-1}}\left(\Omega\right)}^{\frac{2}{p-1}}.
\end{equation}
		
Since $\frac{2}{p-1} < 2$, we have $L^{2}(\Omega) \subset L^{\frac{2}{p-1}}(\Omega)$. Hence, for every $v - v_h$
\begin{equation}\label{Q2}
\left\Vert v-v_h \right\Vert_{L^{\frac{2}{p-1}}(\Omega)}^{\frac{2}{p-1}} \leq C \left\Vert v-v_h \right\Vert_{L^{2}(\Omega)}^{\frac{2}{p-1}}.
\end{equation}
		
By \eqref{Eq03:07_06} and \eqref{Q2}, 	
\begin{equation}\label{Q4}
\left\Vert \left|v - v_h\right|^{\frac{1}{p-1}}\right\Vert_{L^{2}(\Omega)}^2 \leq C \left\Vert v-v_h \right\Vert_{L^{2}(\Omega)}^{\frac{2}{p-1}}.
\end{equation}
	
If $p \geq 3$, then $0 < \frac{2}{p-1} \leq 1$. In this case, take $\mu = \frac{2}{p-1}$ in Theorem \ref{Teo_mu}. Equation $(\ref{Eq01:07_06})$ then becomes
\begin{equation}\label{18/06_Eq11}
\left\Vert \left|v - v_h\right|^{\frac{1}{p-1}}\right\Vert_{L^{2}(\Omega)}^2 \leq  C \left\Vert  v - v_h \right\Vert_{L^{2}(\Omega)}^{\frac{2}{p-1}}.
\end{equation}
	
From \eqref{Q4} and \eqref{18/06_Eq11}, with $2 < p < \infty$, we have
\begin{equation}\label{18/06_Eq12}
\left\Vert \left|v - v_h\right|^{\frac{1}{p-1}}\right\Vert_{L^{2}(\Omega)}^2 \leq C \left\Vert  v - v_h \right\Vert_{L^{2}(\Omega)}^{\frac{2}{p-1}}.
\end{equation}
Now, substituting \eqref{18/06_Eq12} in \eqref{Eq04:07_06} and using \eqref{Ordem_v_p>2}, we obtain the inequality
\begin{equation}\label{Eq05:07_06}
\left\Vert \nabla \theta \right\Vert_{L^{2}(\Omega)}  \leq Ch^{s-1}\left\Vert u \right\Vert_{H^{s}(\Omega)}   + Ch^{\frac{s}{p-1}}\left\Vert v \right\Vert_{H^{s}(\Omega)}^{\frac{1}{p-1}}.
\end{equation}
On the other hand, using the definitions of $\gamma$ and $\theta$ in \eqref{gamma1} and \eqref{theta1}, respectively, and the Minkowski inequality, we get
\begin{equation*}
\left\Vert  \nabla\left(u-u_h\right)\right\Vert_{L^{2}(\Omega)} \leq \left\Vert \nabla \gamma \right\Vert_{L^{2}(\Omega)} + \left\Vert \nabla \theta \right\Vert_{L^{2}(\Omega)}.
\end{equation*}
From Theorem \ref{Teo.Est.Interpolador},
\begin{equation*}
\left\Vert  \nabla\left(u-u_h\right)\right\Vert_{L^{2}(\Omega)} \leq Ch^{s-1}\left\Vert u \right\Vert_{H^{s}(\Omega)} + \left\Vert \nabla \theta \right\Vert_{L^{2}(\Omega)}.
\end{equation*}
Substituting \eqref{Eq05:07_06} in the inequality above, leads to \eqref{Ordem_grad u_p>2}.
	
Assuming $\varphi \in L^{2}(\Omega)$ and $\psi \in H^{2}(\Omega) \cap H_{0}^{1}(\Omega)$ is the solution of $-\Delta \psi = \varphi$ and reasoning as in the proof of Theorem \ref{Teo_Ordem_Conv_1<p<2}, we obtain \eqref{Ordem_u_p>2} which completes the proof.
\end{proof}

\subsection{Discretization by the finite element method}

In this section, our aim is to transform Problem \eqref{Sistema_discreto} into an algebraic system of equations. For this purpose, we consider the following basis for the finite element space
\begin{align*}
V_{h} = \langle\varphi_1, \varphi_2, \cdots, \varphi_n\rangle .
\end{align*}

Then, we can represent
\begin{equation}  \label{u_h, v_h, psi_h, varphi_h}
\begin{array}{llll}
u_h = \displaystyle\sum_{i=1}^{n}u_i\varphi_i, & 
v_h = \displaystyle\sum_{i=1}^{n}v_i\varphi_i, & 
\eta_h = \displaystyle\sum_{i=1}^{n}\eta_i \varphi_i, & 
\psi_h = \displaystyle\sum_{i=1}^{n}\psi_i \varphi_i,
\end{array}
\end{equation}
where $u_i$ and $v_i$, with $1 \leq i \leq n$, are the coefficients to be determined. Substituting \eqref{u_h, v_h, psi_h, varphi_h} in Problem \eqref{Sistema_discreto} and noting that $\psi_i$ and $\eta_i$ are arbitrary, we obtain the system
\begin{equation*}
\left\{
\begin{array}{ll}
\left( \nabla \displaystyle\sum_{i=1}^{n}v_i\varphi_i, \nabla \varphi_j \right) = -\left(f, \varphi_j\right), & j = 1, \cdots, n, 
\vspace{0.1cm} \\
\left( \nabla \displaystyle\sum_{i=1}^{n}u_i\varphi_i, \nabla \varphi_j\right) = - \left( \left|\displaystyle\sum_{i=1}^{n}v_i\varphi_i \right|^{q-2}\displaystyle\sum_{i=1}^{n}v_i\varphi_i, \varphi_j\right), & j = 1, \cdots, n.
\end{array}
\right.
\end{equation*}
Since $u_{i}$ and $v_{i}$, $1\leq i\leq n$, are constant coefficients, we have
\begin{equation}  \label{Sistema6}
\left\{
\begin{array}{ll}
\displaystyle\sum_{i=1}^{n}v_i\left( \nabla \varphi_i, \nabla\varphi_j \right) = -\left( f, \varphi_j \right), & j = 1, \cdots, n, 
\vspace{0.1cm} \\
\displaystyle\sum_{i=1}^{n}u_i\left( \nabla \varphi_i, \nabla \varphi_j \right) = - \left( \left|\displaystyle\sum_{i=1}^{n}v_i\varphi_i \right|^{q-2}\displaystyle\sum_{i=1}^{n}v_i\varphi_i, \varphi_j\right), & j = 1, \cdots, n.
\end{array}
\right.
\end{equation}
Notice that
\begin{equation*}
\left( \left|\displaystyle\sum_{i=1}^{n}v_i\varphi_i \right|^{q-2}\displaystyle\sum_{i=1}^{n}v_i\varphi_i, \varphi_j\right) = \displaystyle\sum_{i=1}^{n}v_i\int_{\Omega} \left|\displaystyle\sum_{k=1}^{n}v_k\varphi_k \right|^{q-2}\varphi_i\varphi_j dx.
\end{equation*}
We denote the elements of vectors $\vec{u} \in \mathbb{R}^{n \times 1}$ and $\vec{v} \in \mathbb{R}^{n \times 1}$ by $u_i$ and $v_i$, respectively, while the elements of matrices $M(\vec{v}) \in \mathbb{R}^{n \times n}$ and $K \in \mathbb{R}^{n \times n}$ are denoted by
\begin{align*}
& M_{ij} = \displaystyle\int_{\Omega} \left|\displaystyle\sum_{k=1}^{n}v_k \varphi_k \right|^{q-2}\varphi_i\varphi_j dx, 
\\
& K_{ij} = \left( \nabla\varphi_i, \nabla \varphi_j \right),
\end{align*}
respectively. The elements of vector $\vec{f} \in \mathbb{R}^{n \times 1}$ are denoted by $f_i = -\left( f, \varphi_i \right)$.  Thus, we obtain the following system associated with Problem \eqref{Sistema6}
\begin{equation}  \label{System}
\left\{
\begin{array}{l}
K\vec{v} = \vec{f}, 
\\
K\vec{u} = - M\left(\vec{v}\right)\vec{v}.
\end{array}
\right.
\end{equation}

Theorems \ref{Teo_existencia_uh_vh}, \ref{Teo_unicidade_uh_vh} and \ref{Teo_estabilidade_vh_uh} guarantee the existence, uniqueness and regularity of solutions $\vec{u}$ and $\vec{v}$ of System \eqref{System}. In addition, since  matrix $K$ is formed by the inner product of $\nabla \varphi_i$ and $\nabla \varphi_j$, and the functions $\varphi_i$ and $\varphi_j$ form the finite element space $V_h$, we have that $K$ is tridiagonal, symmetric and positive definite. Therefore, \eqref{System} is a consistent and independent linear system.

\section{Numerical results}\label{sec5}

In this section, we present the results of a Matlab implementation of the theory. First, we validate the code and analyze the order of convergence.

\subsection{Example 1}

Suppose that $\Omega = [0,1]$, the domain mesh is uniform, the finite element space is formed by linear basis functions, $f(x) = 1$ and $p=1.5$. Then,  the exact solutions are
\begin{align*}
u(x) = \dfrac{x - x^4(2x^2 - 6x + 5)}{240}
\end{align*}
\noindent and
\begin{align*}
v(x) = \left|\Delta u(x) \right|^{p-2}\Delta u(x) = \dfrac{x \left(x-1\right)}{2}.
\end{align*}

\begin{figure*}[ht!]
\centering
\begin{subfigure}
[$u(x)$ and $u_h(x)$.]{\label{Ex1_u}\includegraphics[height=7cm]{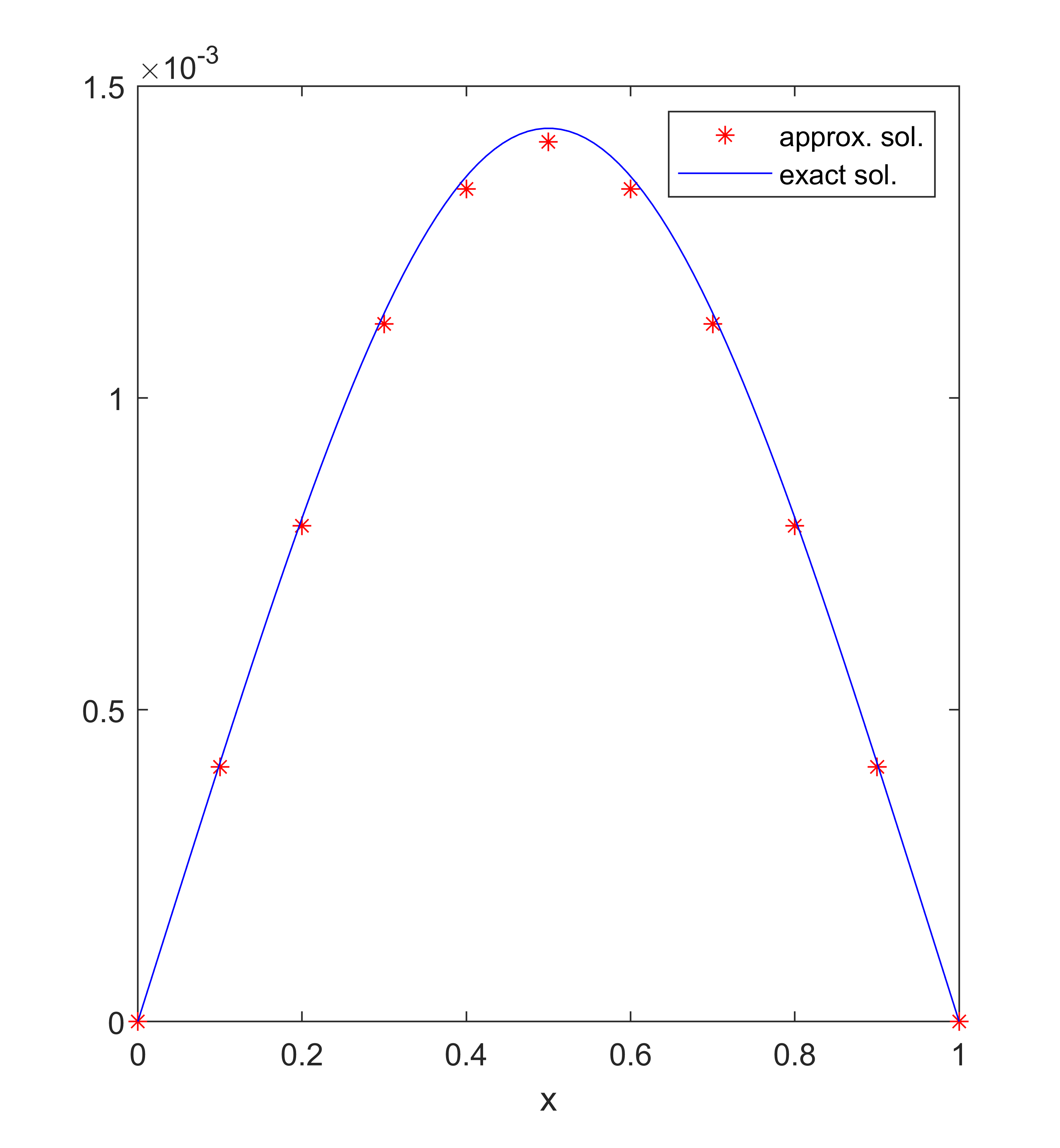}}\end{subfigure}
\begin{subfigure}
[$v(x)$ and $v_h(x)$.]{\label{Ex1_v}\includegraphics[height=7cm]{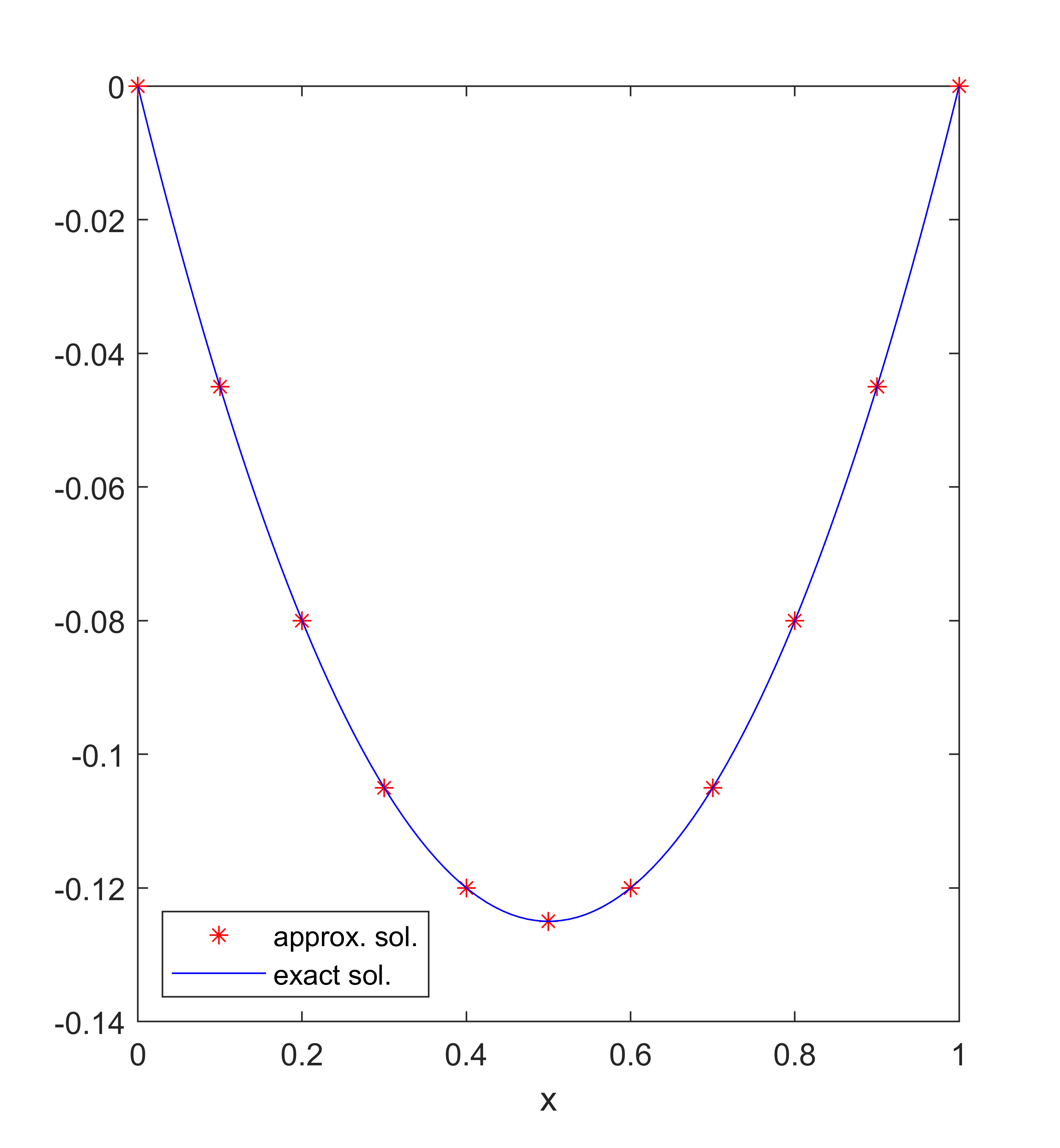}}
\end{subfigure}
\caption{Exact solutions $u(x)$ and $v(x)$ and approximate solutions $u_h(x)$ and $v_h(x)$.}
\end{figure*}

Figures \ref{Ex1_u} and \ref{Ex1_v} show the approximate solutions $u_h(x)$ and $v_h(x)$ compared to the exact solutions $u(x)$ and $v(x)$, respectively. In this case, we use $n=10$ finite elements.

\begin{figure}[H]
\centering
\includegraphics[height=7cm]{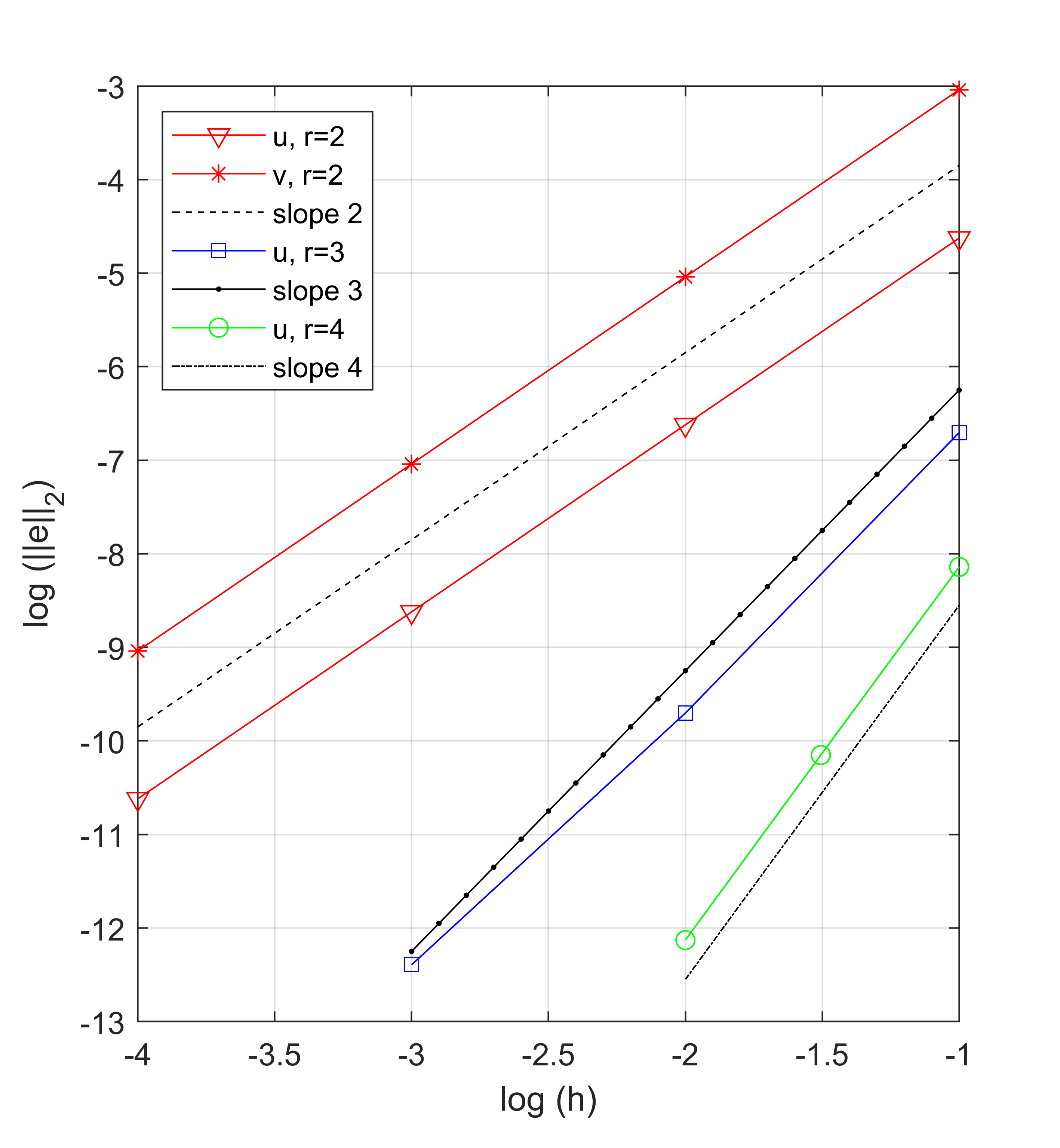}
\caption{Convergence order with $r=2, 3, 4$ and $p=1.5$.}
\label{Order_p=1.5}
\end{figure}

In Figure \ref{Order_p=1.5}, we show the graphs of the order of convergence for $p=1.5$. We considered $n=10, 100, 1000, 10000$ and, for each case, we calculated the error in the $L^2(\Omega)$ norm.

We note that function $v(x)$ is a polynomial of degree $2$ and so, for $r \geq 3$, the approximate solution is identical to the exact solution. Hence the errors depend only on the method used for solving the system, which means that there is no order of convergence. Therefore, in Figure \ref{Order_p=1.5}, the convergence order graph for $v_h(x)$ is not displayed when using base functions with degree $2$ and $3$, respectively, $r=3$ and $r=4$.

Next, we still consider $f = 1$, but different values of $p$. For each value of $p$, we calculate the exact solution. At this stage, the approximate solutions $u_h(x)$ and $v_h(x)$ are very close to the exact solutions $u(x)$ and $v(x)$, so we choose not to display the graphs that perform the comparison between the approximate and exact solutions.

\begin{figure}[H]
\centering
\begin{subfigure}[$p=1.1$.]{\label{Order_p=1.1}
\includegraphics[height=7cm]{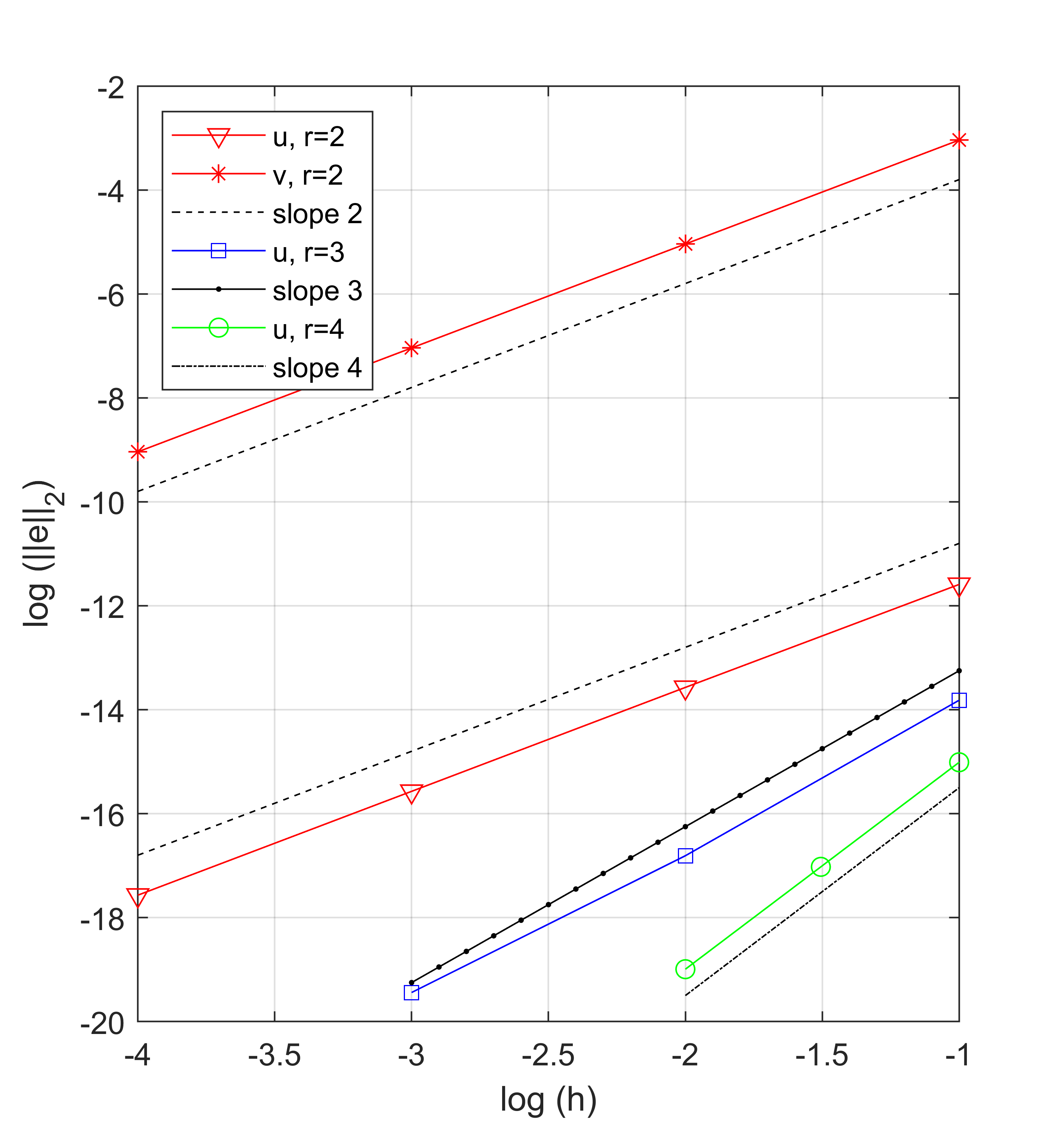}}  
\end{subfigure}
\begin{subfigure}
[$p=2$.]{\label{Order_p=2}\includegraphics[height=7cm]{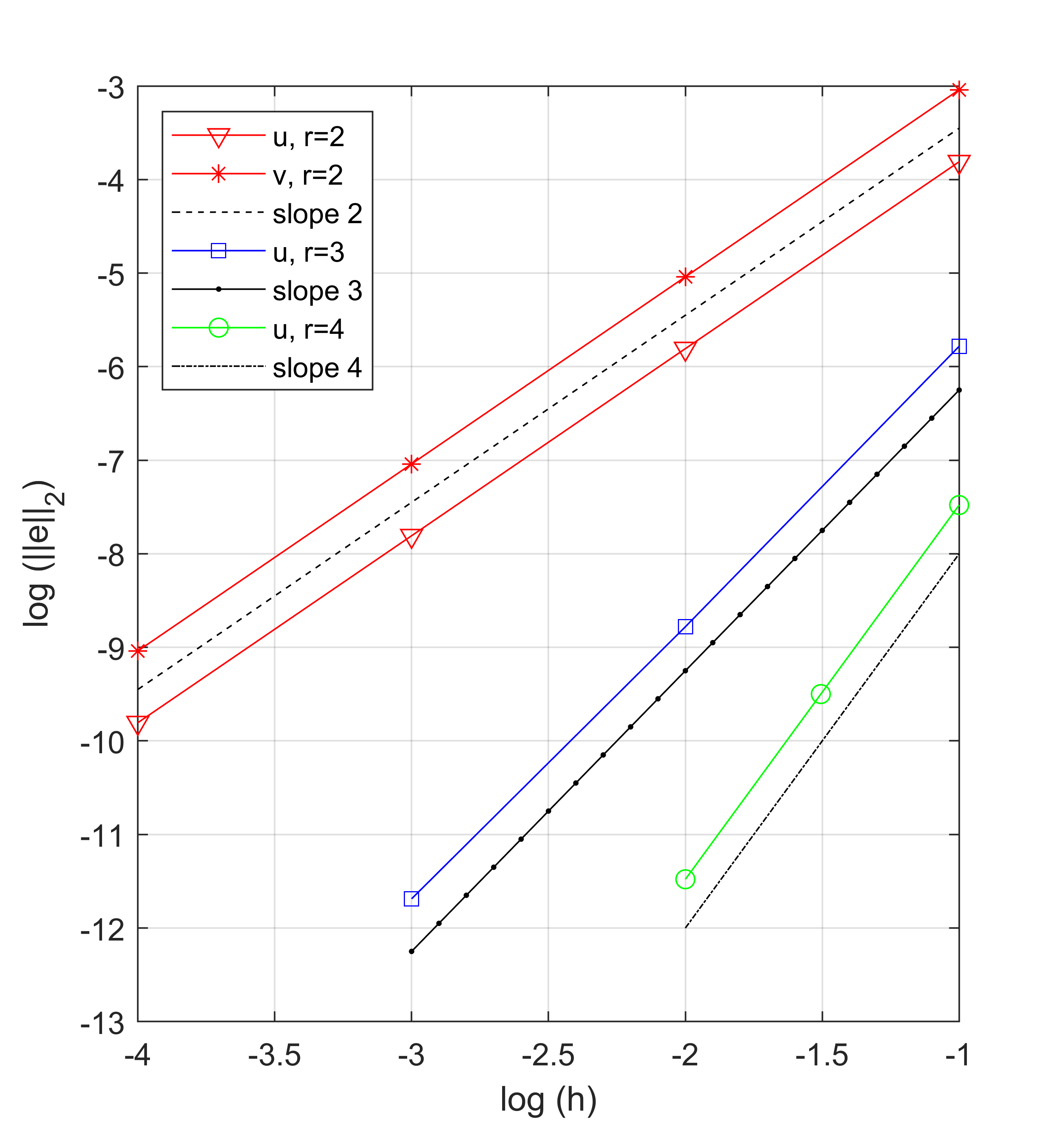}}
\end{subfigure}
\caption{Convergence order with $r=2, 3, 4$.}
\label{Order_p=1.1 and 2}
\end{figure}

We emphasize that, in Figures \ref{Order_p=1.5}, \ref{Order_p=1.1} and \ref{Order_p=2}, the convergence order graphs are in agreement with the results of Theorem \ref{Teo_Ordem_Conv_1<p<2}, that is, the computer simulations are consistent with the analytical study. Furthermore, estimates \eqref{Ordem_v_1<p<2} and \eqref{Ordem_u_1<p<2} do not depend on $p$ and, for $\Vert u - u_h \Vert_{L^2\left(\Omega\right)}$ and $\Vert v - v_h \Vert_{L^2\left(\Omega\right)}$, the convergence order is $O(h^r)$ for polynomials of degree $r-1$, with $r \geq 2$.

\subsection{Example 2}

In this example, suppose that $\Omega = [0,1]$ and consider then the domain mesh is uniform, the finite element space is formed by linear basis functions, $f(x) = -\frac{5x^4}{6} + \frac{2x^2}{3} - \frac{1}{18}$ and $p=3$. In this case, the exact solutions are
\begin{align*}
u(x) = \dfrac{x^5}{120} - \dfrac{x^3}{36} + \dfrac{7x}{360}
\end{align*}
\noindent and
\begin{align*}
v(x) = \dfrac{x^2\left(x^4 -2x^2 + 1\right)}{36}.
\end{align*}

\begin{figure}[H]
\centering
\begin{subfigure}[$u(x)$ and $u_h(x)$.]{\label{Ex2_u}\includegraphics[height=7cm]{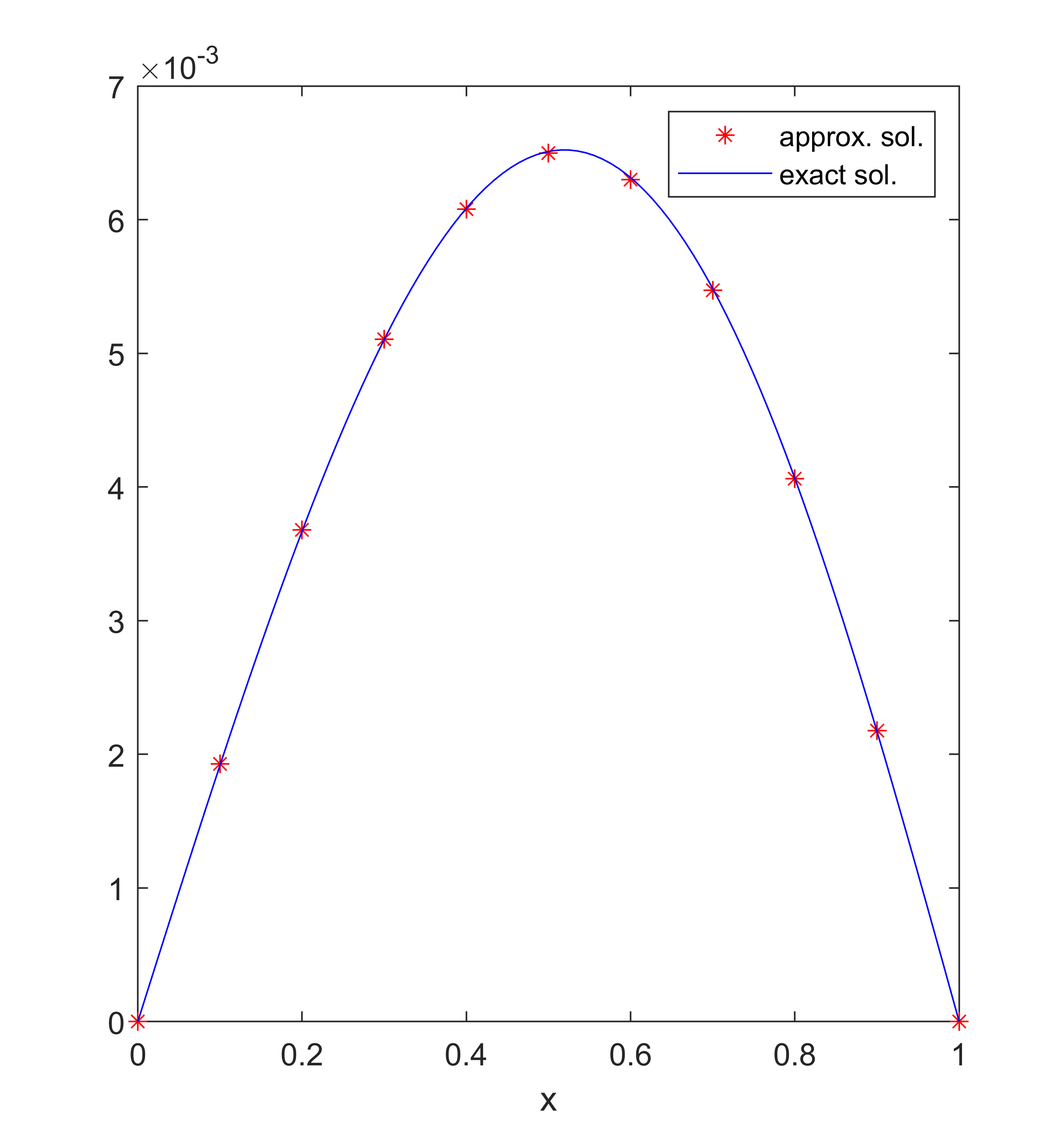}} 
\end{subfigure}
\begin{subfigure}[$v(x)$ and $v_h(x)$.]{\label{Ex2_v}\includegraphics[height=7cm]{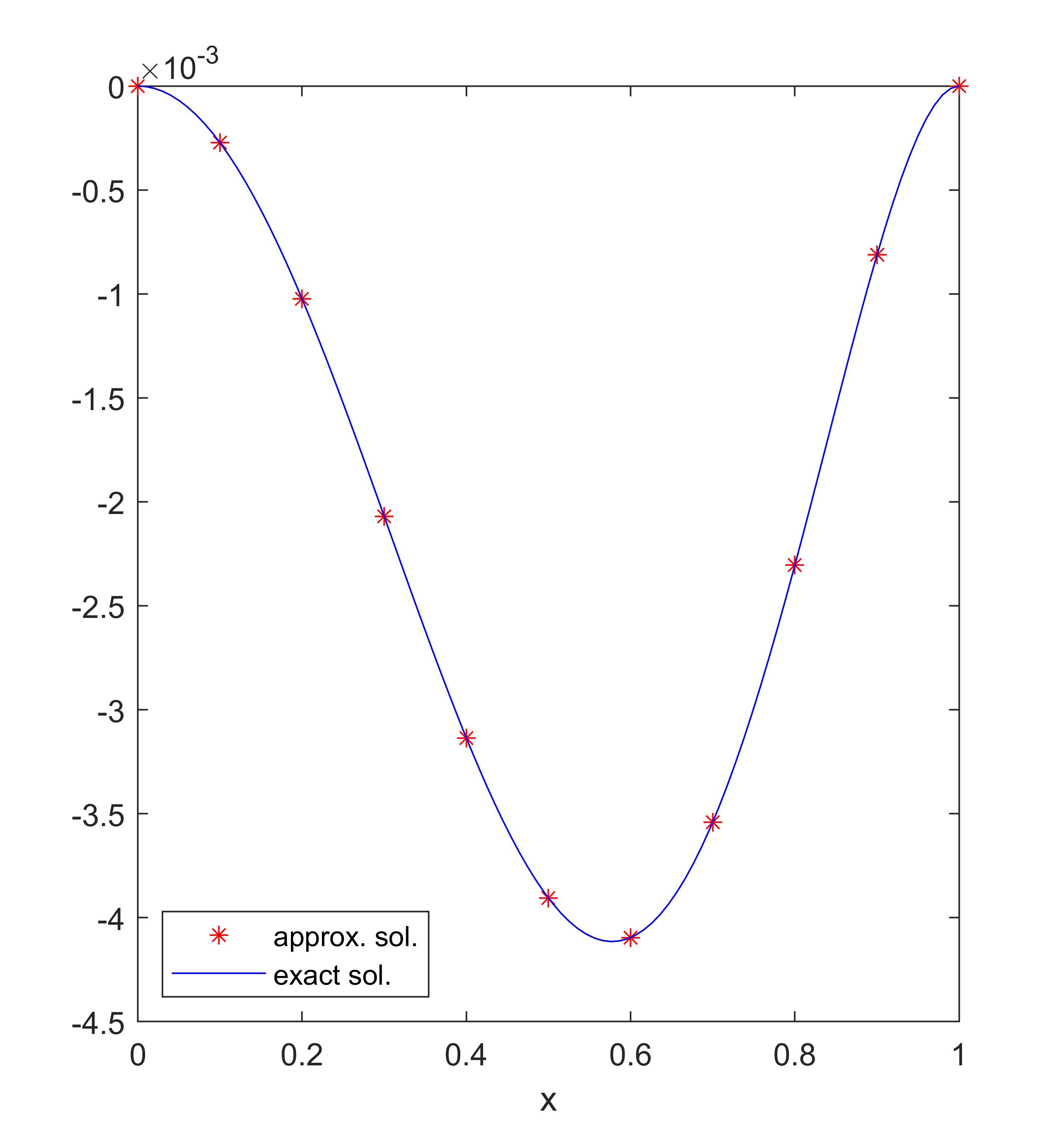}}
\end{subfigure}
\caption{Exact solutions $u(x)$ and $v(x)$ and approximate solutions $u_h(x)$ and $v_h(x)$.}
\end{figure}

Figures \ref{Ex2_u} and \ref{Ex2_v} compare the approximate solutions $u_h(x)$ and $v_h(x)$ to the exact solutions $u(x)$ and $v(x)$, respectively. In this case, $n=10$ finite elements are also used.

\begin{figure}[H]
\centering
\includegraphics[height=7cm]{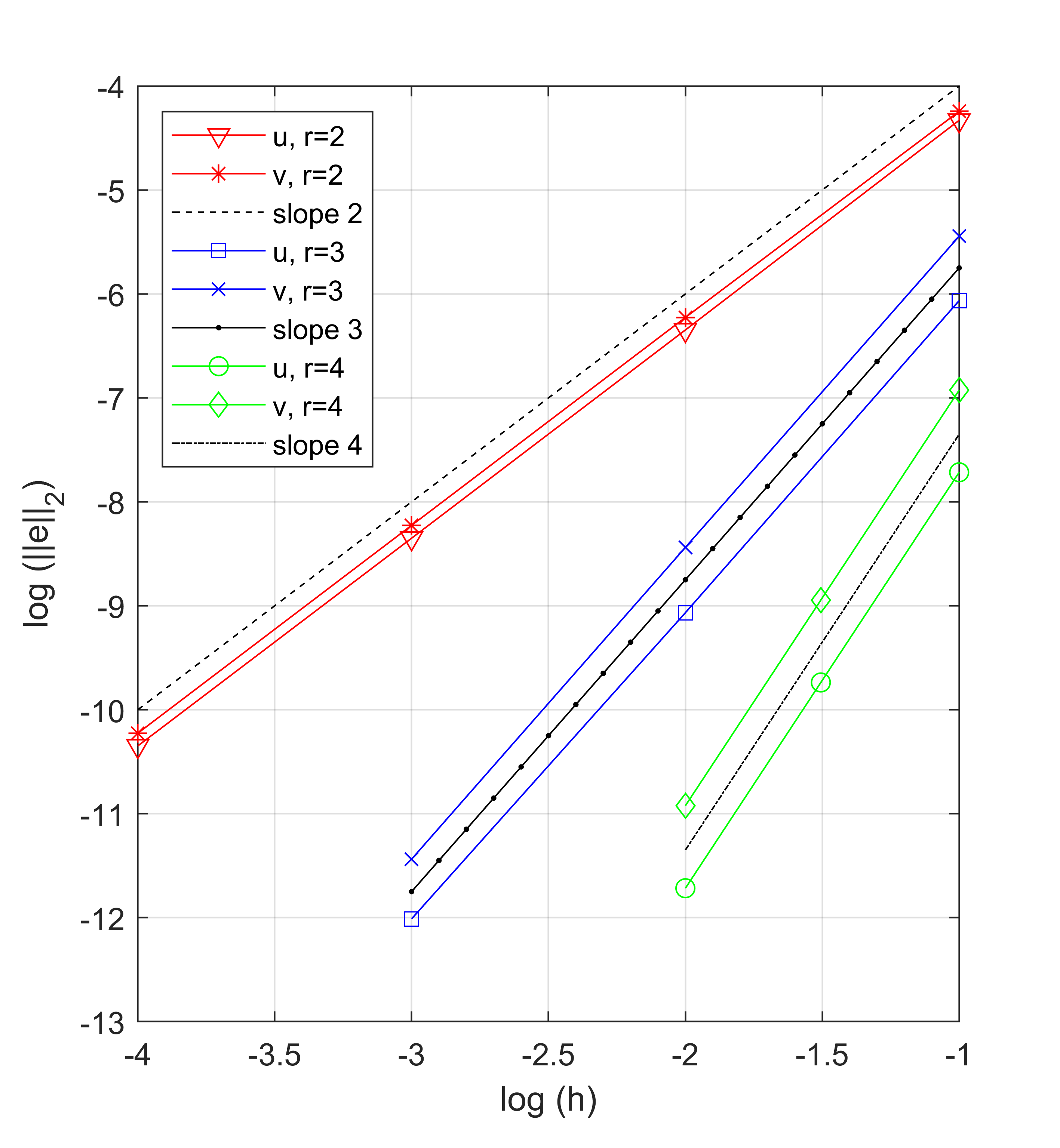}
\caption{Convergence order with $r=2, 3, 4$ and $p=3$.}
\label{Order_p=3}
\end{figure}

In Figure \ref{Order_p=3}, we show the graphs of the order of convergence for $p=3$. Again, we considered $n=10, 100, 1000, 10000$ and, for each one, we calculated the error in the $L^2(\Omega)$ norm.

Now, we change the values of $p$ and obtain the function $f(x)$ given by
\begin{equation*}
f(x) = x\left(p-1\right) \left(\dfrac{x}{6}-\dfrac{x ^3}{6}\right)^{p-2}-\left(p-1\right)\left(p-2\right)\left(\frac{x}{6}-\frac{x^3}{6}\right)^{p-3}\left(\frac{x^2}{2}-\frac{1}{6}\right)^2.
\end{equation*}

In this case, the approximate solutions are very close to the exact solutions, so we omit the graphs that compare the approximate and exact solutions.

\begin{figure}[H]
\centering
\begin{subfigure}[$p=4$.]{\label{Order_p=4}\includegraphics[height=7cm]{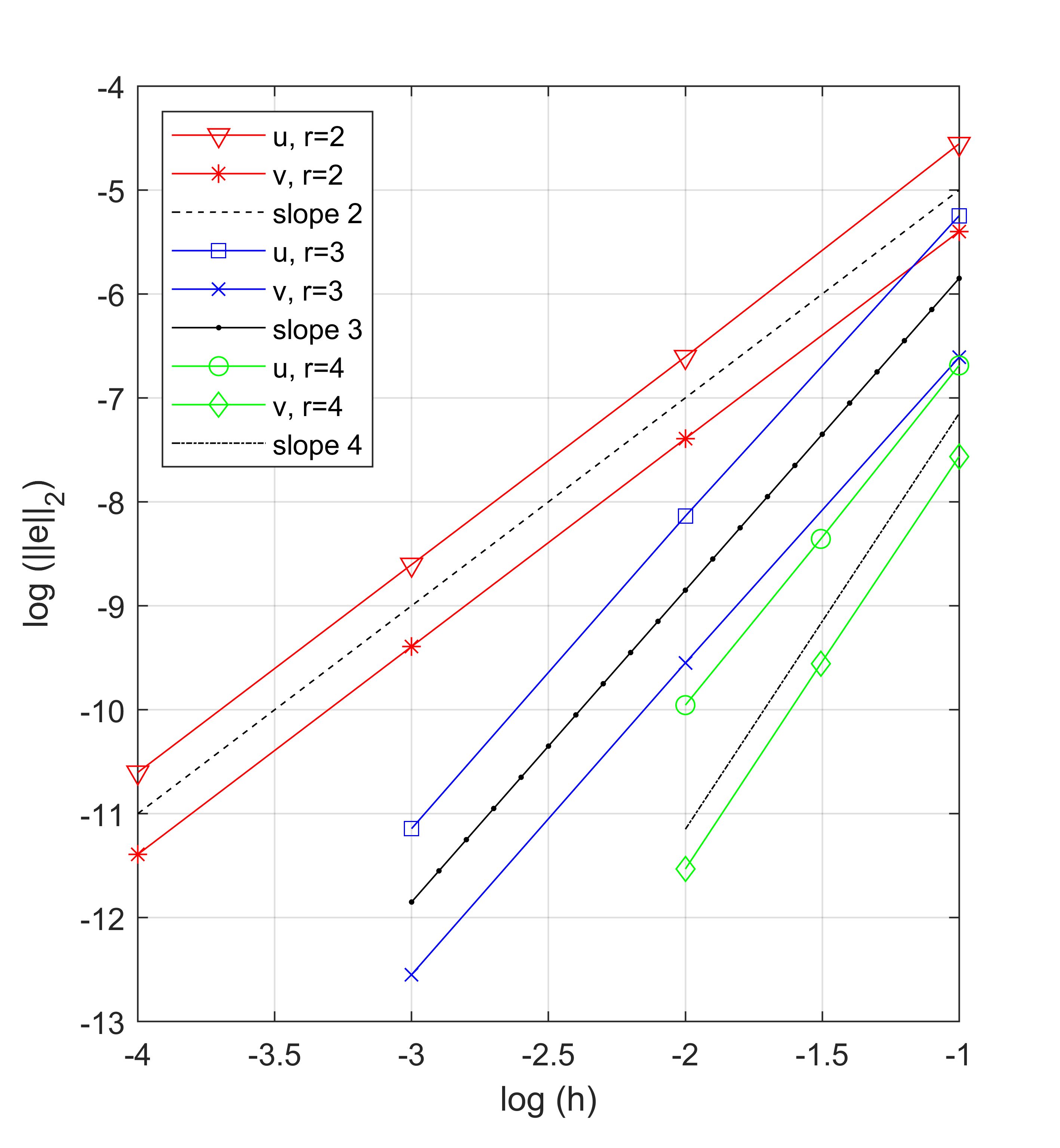}}
\end{subfigure}
\begin{subfigure}[$p=5$.]{\label{Order_p=5}\includegraphics[height=7cm]{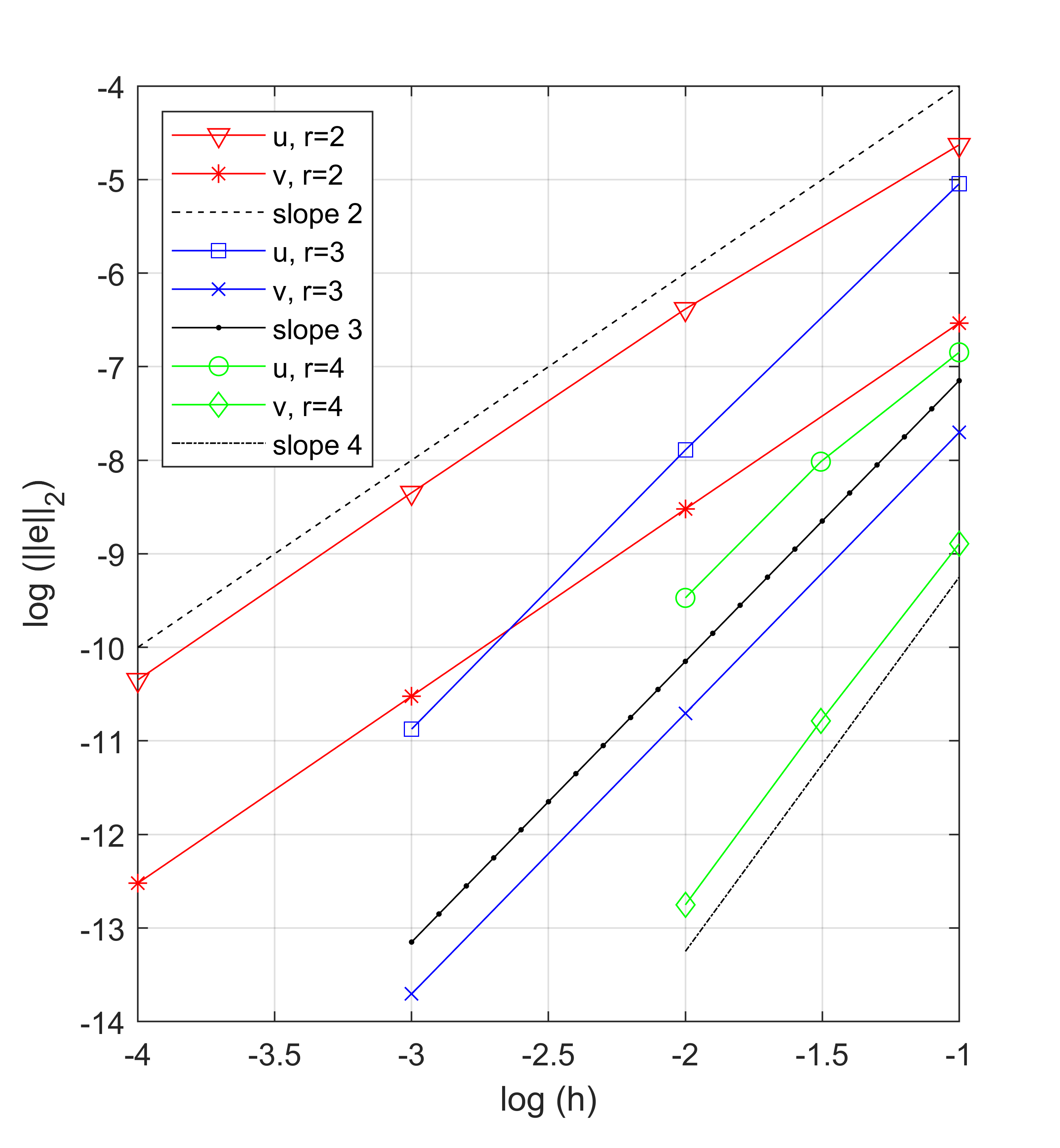}}
\end{subfigure}
\caption{Convergence order with $r=2, 3, 4$.}
\label{Order_p=4 and 5}
\end{figure}

\begin{figure}[H]
\centering
\begin{subfigure}[$p=10$.]{\label{Order_p=10_u}\includegraphics[height=7cm]{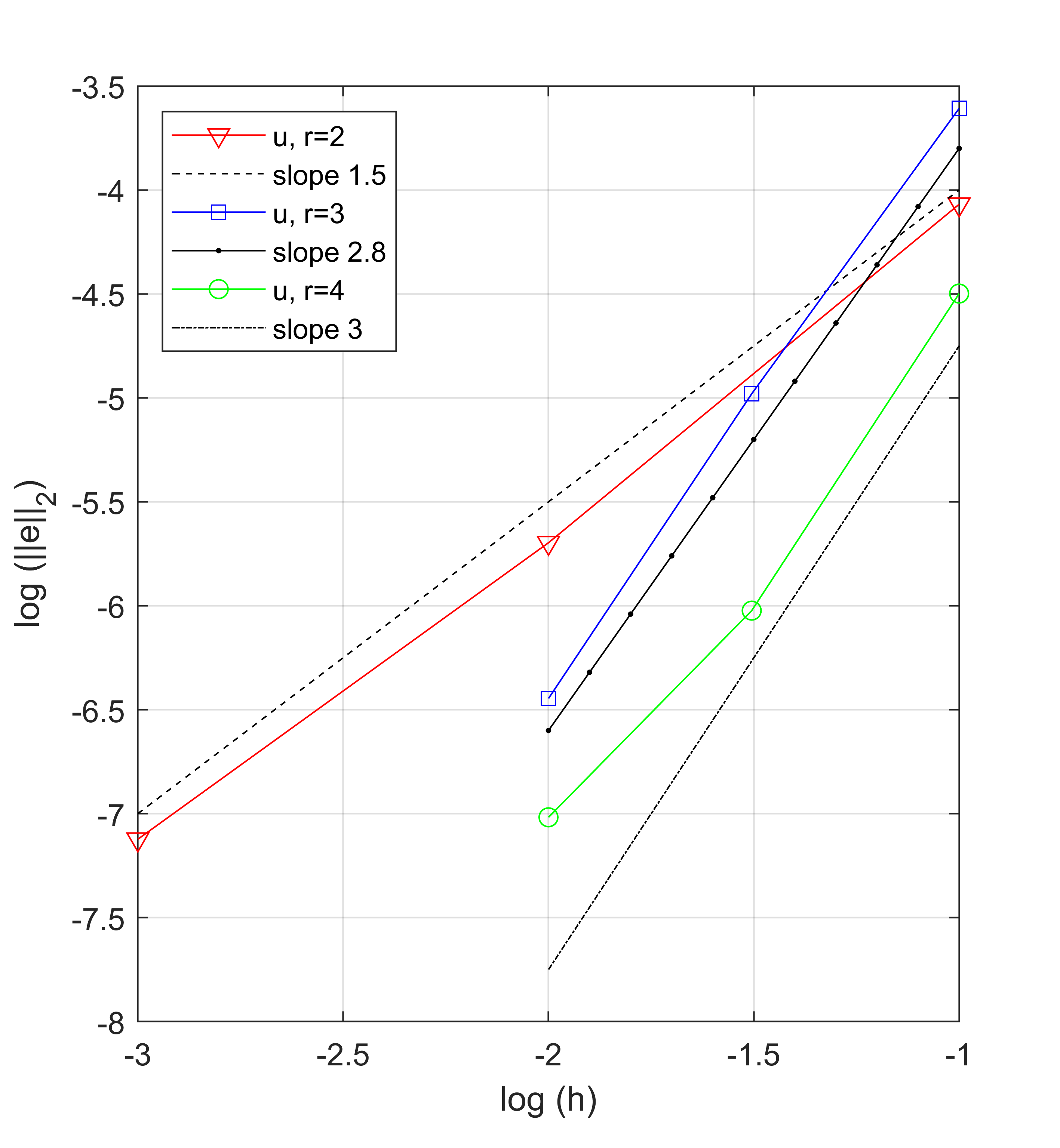}}
\end{subfigure}
\begin{subfigure}[$p=10$.]{\label{Order_p=10_v}\includegraphics[height=7cm]{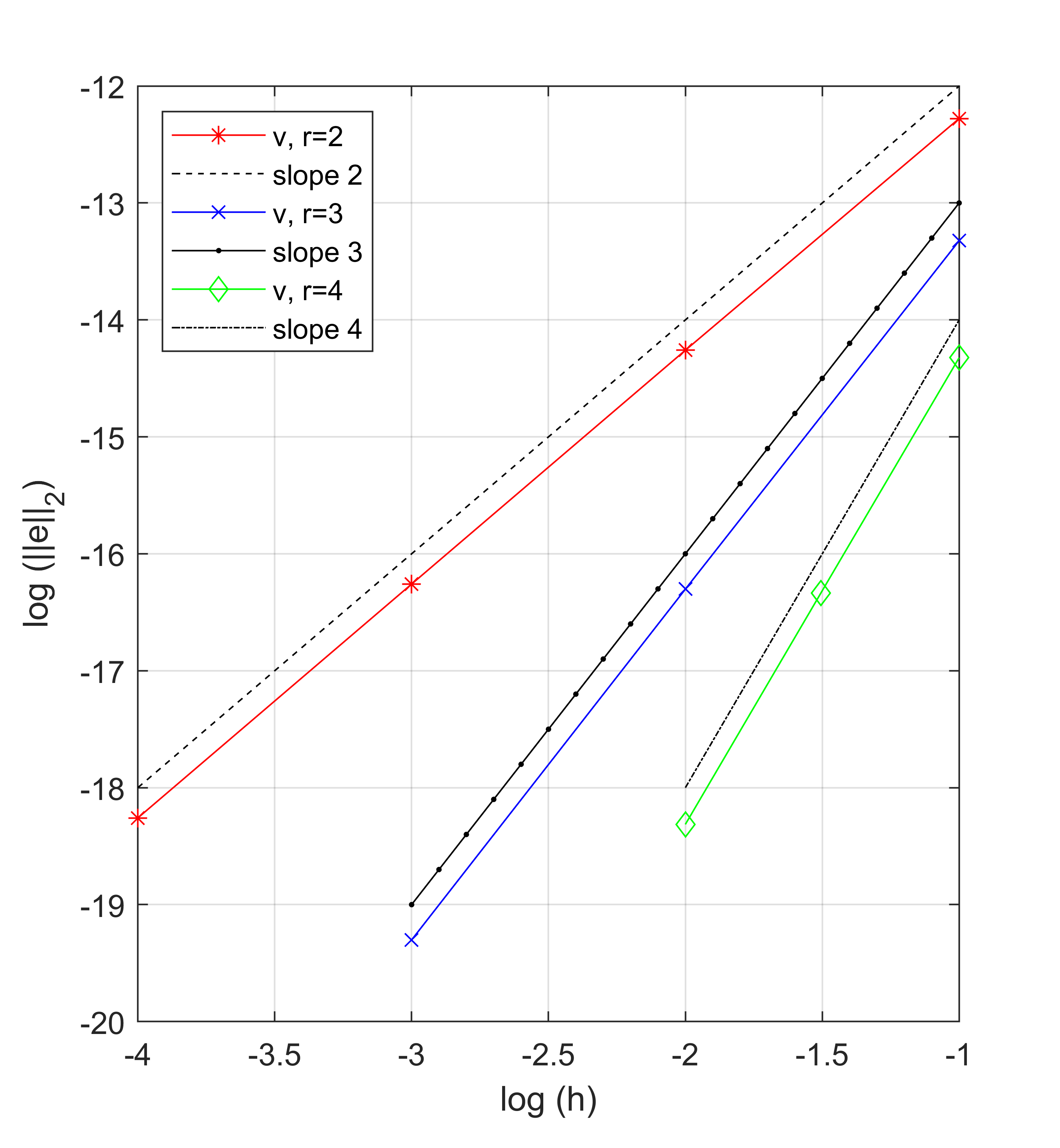}}
\end{subfigure}
\caption{Convergence order with $r=2, 3, 4$.}
\label{Order_p=10}
\end{figure}

\begin{figure}[H]
\centering
\begin{subfigure}[$p=25$.]{\label{Order_p=25_u}\includegraphics[height=7cm]{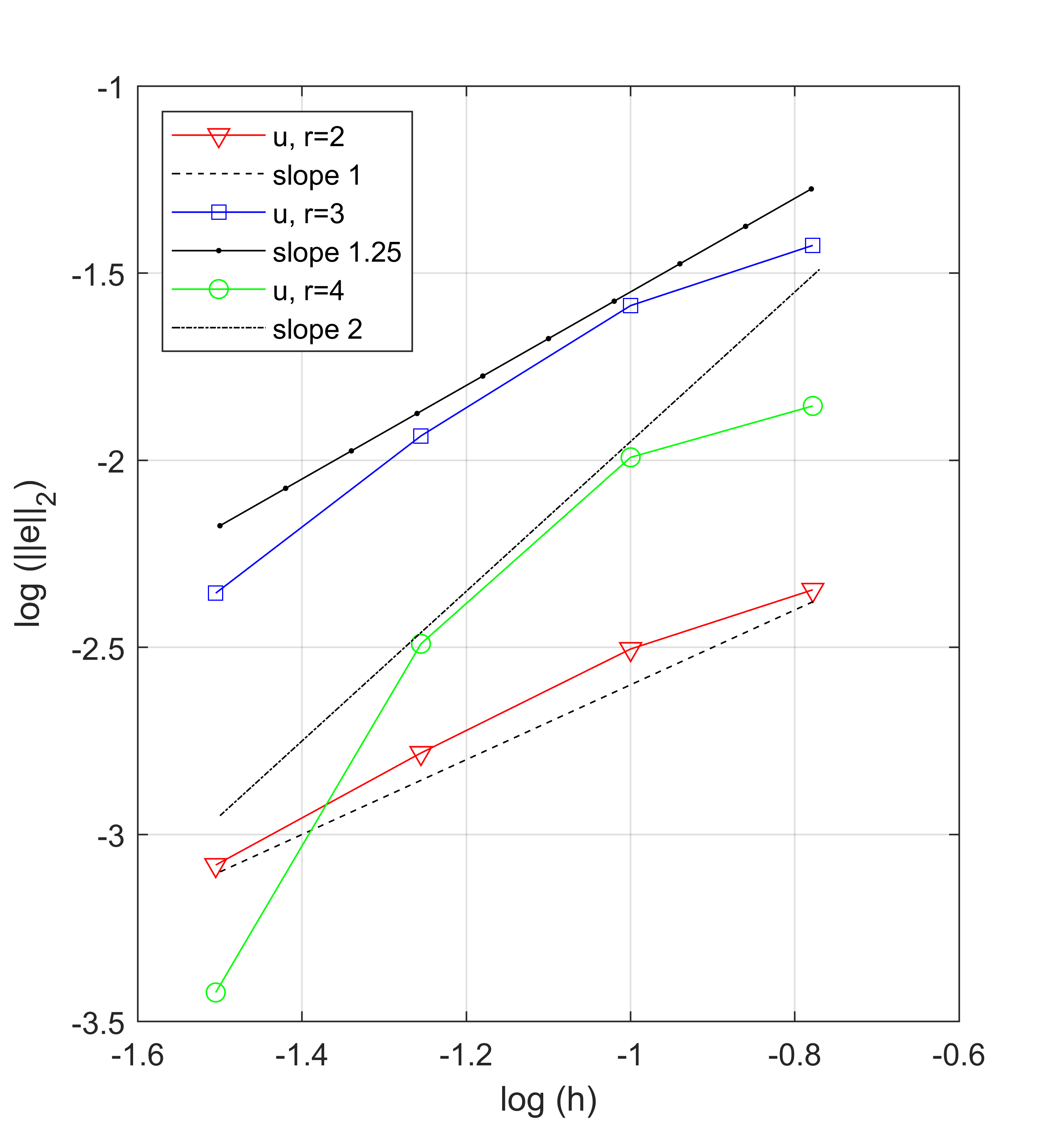}}
\end{subfigure}
\begin{subfigure}[$p=25$.]{\label{Order_p=25_v}\includegraphics[height=7cm]{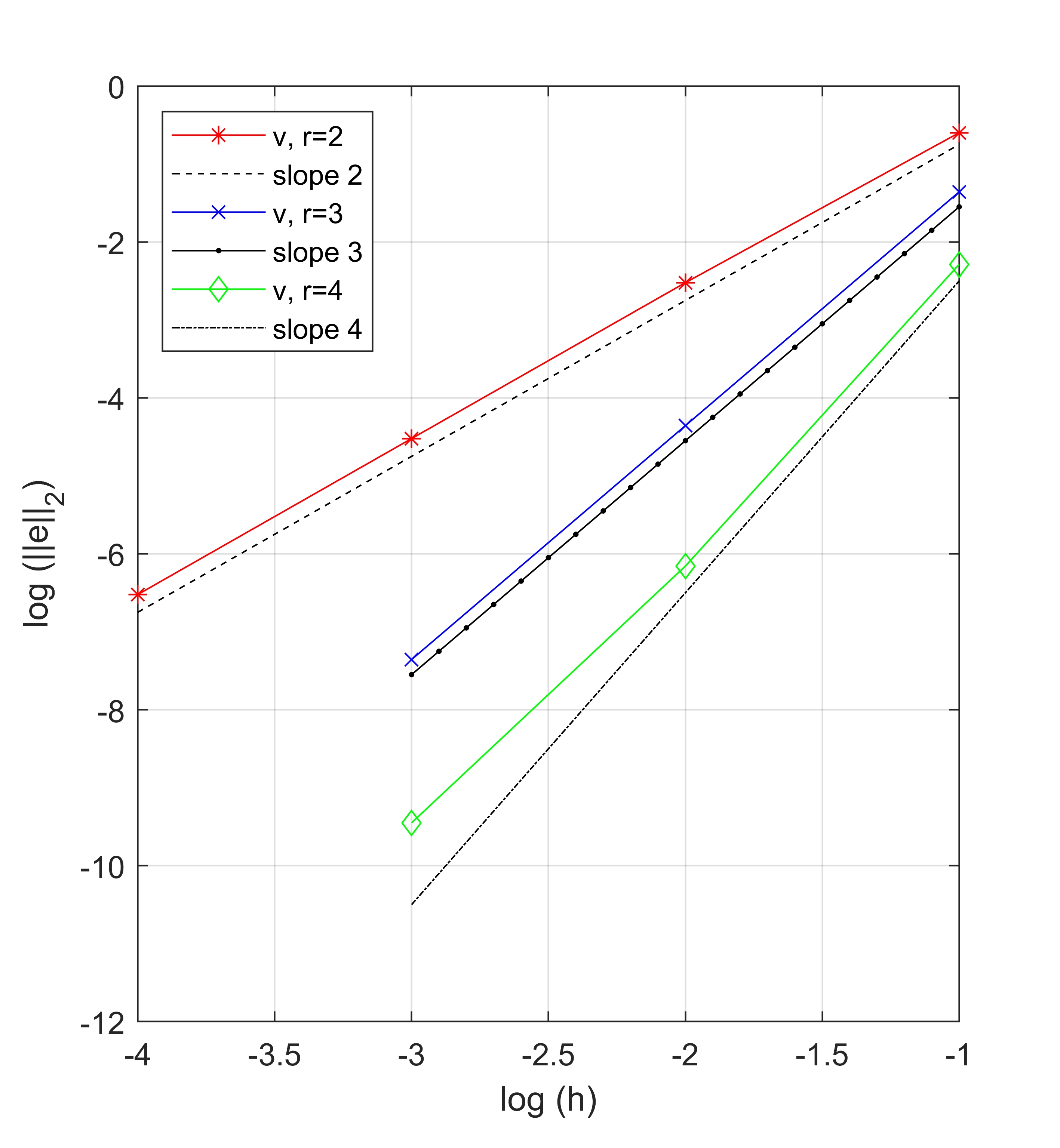}}
\end{subfigure}
\caption{Convergence order with $r=2, 3, 4$.}
\label{Order_p=25}
\end{figure}

Let us recall that, estimate \eqref{Ordem_v_p>2} in Theorem \ref{Teo_Ordem_Conv_p>2} does not depend on $p$, so in Figures \ref{Order_p=3}, \ref{Order_p=4}, \ref{Order_p=5}, \ref{Order_p=10_v} and \ref{Order_p=25_v}, we have that, for $\Vert v - v_h \Vert_{L^2\left(\Omega\right)}$, the convergence order is $O(h^r)$ for polynomials of degree $r-1$, with $r \geq 2$.

In Figure \ref{Order_p=3}, for $\Vert u - u_h\Vert_{L^{2}\left(\Omega\right)}$, we obtain that the convergence order is $O(h^r)$ when we use polynomials of degree $r-1$, with $r \geq 2$. The analysis is identical for Figures \ref{Order_p=4} and \ref{Order_p=5}, but, for $r=4$, we note that the convergence order is less than $O(h^4)$. For Figures \ref{Order_p=10_u} and \ref{Order_p=25_u}, the convergence order exists but is low, because in Theorem \ref{Teo_Ordem_Conv_p>2} the estimate \eqref{Ordem_u_p>2} depends on $p$ and when $p \rightarrow +\infty$ the order of convergence decreases.

\section{Conclusions}\label{sec6}

We study a nonlinear beam equation with the $p$-biharmonic operator. Rewriting Problem \eqref{Problema} as a system of adequate differential equations, we prove the existence, uniqueness and regularity of the solution. Moreover, we prove the existence, uniqueness and stability of the discrete solution. We established sufficient conditions on the data to obtain optimal convergence rates for $1 < p < 2$ and some finite element solutions with piecewise polynomial of arbitrary degree basis functions in space. Finally, we implement the method in Matlab software for $ N=1$ and perform some simulations that confirm the theory.

An interesting idea is to generalize Problem \eqref{Problema} by using the variable exponent $p(x)$. In this case, it will be necessary to use the definitions and properties of the Lebesgue and Sobolev spaces with variable exponents, so the study of the solution will not be trivial.

\section*{Acknowledgments}

This work was partially supported by the research projects: FEDER through the Programa Operacional Factores de Competitividade, FCT - Fundação para a Ciência e a Tecnologia [Grant Number UIDB/00212/2020].

The fourth author was supported by FCT - Fundação para a Ciência e a Tecnologia, through Centro de Matemática e Aplicações - Universidade da Beira Interior, under the Grant Number  UI/BD/150794/2020, and also supported by MCTES, FSE and UE.


\end{document}